\documentclass[AMA,STIX1COL]{WileyNJD-v2}
\pdfoutput=1
\usepackage{moreverb}
\usepackage{epsfig} 
\usepackage{graphicx}
\usepackage[ansinew]{inputenc}
\usepackage{epstopdf}
\usepackage{psfrag}
\usepackage{amsmath}
\usepackage{amssymb}
\usepackage{float}
\usepackage{multirow}
\usepackage{lineno}
\usepackage{setspace}
\usepackage{flushend}
\usepackage{array}
\usepackage{tabu}
\usepackage{bm}
\usepackage{color}
\usepackage{url}
\newcolumntype{L}[1]{>{\raggedright\let\newline\\\arraybackslash\hspace{0pt}}m{#1}}
\newcolumntype{C}[1]{>{\centering\let\newline\\\arraybackslash\hspace{0pt}}m{#1}}
\newcolumntype{R}[1]{>{\raggedleft\let\newline\\\arraybackslash\hspace{0pt}}m{#1}}
\usepackage[figuresright]{rotating}
\usepackage{subfigure}
\usepackage{amsfonts}
\usepackage{svgcolor}
\usepackage{styles/ieeetrantools}
\usepackage{makecell}
\usepackage{booktabs}
\usepackage{fixmath}
\usepackage{scalerel}
\usepackage{eurosym}
\usepackage{gensymb}

\usepackage{xcolor}
\definecolor{verde}{rgb}{0.0, 0.5, 0.0}
\definecolor{rojo}{rgb}{0.7, 0, 0.0}
\definecolor{grisClaro}{rgb}{0.65, 0.65, 0.65}

\newcommand\BibTeX{{\rmfamily B\kern-.05em \textsc{i\kern-.025em b}\kern-.08em
T\kern-.1667em\lower.7ex\hbox{E}\kern-.125emX}}

\articletype{Research Article}

\received{<day> <Month>, <year>}
\revised{<day> <Month>, <year>}
\accepted{<day> <Month>, <year>}

\begin{document}

\title{MINLP-based hybrid strategy for operating mode selection of TES-backed-up refrigeration systems}

\author[1]{Guillermo Bejarano*}

\author[2]{David Rodr{\'\i}guez}

\author[3]{Jo{\~a}o M. Lemos}

\author[2]{Manuel Vargas}

\author[2]{Manuel G. Ortega}

\authormark{GUILLERMO BEJARANO \textsc{et al.}}

\address[1]{Universidad Loyola Andaluc{\'\i}a, Escuela T{\'e}cnica Superior de Ingenier{\'\i}a, Departamento de Ingenier{\'\i}a, Sevilla, Espa\~na}

\address[2]{Departamento de Ingenier{\'\i}a de Sistemas y Autom{\'a}tica, Escuela T{\'e}cnica Superior de Ingenier{\'\i}a, Universidad de Sevilla (Espa\~na)}

\address[3]{INESC-ID, Instituto Superior T{\'e}cnico, Universidade de Lisboa (Portugal)}

\corres{*Guillermo Bejarano, Universidad Loyola Andaluc{\'\i}a, Escuela T{\'e}cnica Superior de Ingenier{\'\i}a, Departamento de Ingenier{\'\i}a, Campus de Palmas Altas, 41014, Sevilla, Espa\~na. \email{gbejarano@uloyola.es}}

\presentaddress{Universidad Loyola Andaluc{\'\i}a, Escuela T{\'e}cnica Superior de Ingenier{\'\i}a, Departamento de Ingenier{\'\i}a, Campus de Palmas Altas, 41014. Sevilla, Espa\~na.}

\abstract[Abstract]{This brief deals with the satisfaction of the daily cooling demand by a hybrid system that consists of a vapour-compression refrigeration cycle and a thermal energy storage (TES) unit, based on phase change materials. The addition of the TES tank to the original refrigeration plant allows to schedule the cooling production regardless of the instantaneous demand, given that the TES tank can store cold energy and release it whenever deemed appropriate. The scheduling problem is posed as an optimization problem based on mixed-integer non-linear programming (MINLP), since it includes both discrete and continuous variables. The latter corresponds to the references on the main cooling powers involved in the problem (cooling production at the evaporator and TES charging/discharging), whereas the discrete variables define the operating mode scheduling. Therefore, in addition to the hybrid features of the physical plant, a hybrid optimal control strategy is also proposed. A receding horizon approach is applied, similar to model predictive control (MPC) strategies, while economic criteria are imposed in the objective function, as well as feasibility issues. The TES state estimation is also addressed, since its instantaneous \emph{charge ratio} is not measurable. The proposed strategy is applied in simulation to a challenging cooling demand profile and the main advantages of the MINLP-based strategy over a non-linear MPC-based scheduling strategy previously developed are highlighted, regarding operating cost, ease of tuning, and ability to adapt to cooling demand variations.}

\keywords{Refrigeration system; Thermal energy storage; Phase change materials; Mixed-integer non-linear programming; Scheduling}

\jnlcitation{\cname{%
		\author{Bejarano G}, 
		\author{Rodr{\'\i}guez D}, 
		\author{Lemos JM}, 
		\author{Vargas M}, and
		\author{Ortega MG}} (\cyear{2019}), 
	\ctitle{MINLP-based hybrid strategy for operating mode selection of TES-backed-up refrigeration systems},
	\cjournal{Int. J. Robust Nonlinear Control}, \cvol{2019;00:1--6}.}

\maketitle

\footnotetext{\textbf{Abbreviations:} MINLP, mixed-integer non-linear programming; TES, thermal energy storage; MPC, model predictive control.}

\footnotetext{This is the peer reviewed version of the following article: Bejarano, G., Rodríguez, D., Lemos, J. M., Vargas, M., Ortega, M. G. (2020). MINLP-based hybrid strategy for operating mode selection of TES-backed-up refrigeration systems. International Journal of Robust and Nonlinear Control, 30, 6091-6111, which has been published in final form at \url{https://onlinelibrary.wiley.com/doi/10.1002/rnc.4674}. This article may be used for non-commercial purposes in accordance with Wiley Terms and Conditions for Use of Self-Archived Versions. This article may not be enhanced, enriched or otherwise transformed into a derivative work, without express permission from Wiley or by statutory rights under applicable legislation. Copyright notices must not be removed, obscured or modified. The article must be linked to Wiley’s version of record on Wiley Online Library and any embedding, framing or otherwise making available the article or pages thereof by third parties from platforms, services and websites other than Wiley Online Library must be prohibited.}

%%%%%%%%%%%%%%%%%%%%%%%%%%%%%%%%%%%%%%%%%%%%%%%%%%%%%%%%%%%%%%%%%%%%%%%%%%%%%%%%%%%%%

\section{Introduction} \label{secIntroduction}

Refrigeration cycles based on vapour compression constitute the worldwide leading technology for cooling issues, \emph{i.e.} air conditioning, medium-temperature refrigeration, and freezing. Very different areas demand controlling room temperature, for instance for human comfort, food storage and transportation, industrial processes, etc., where a wide power range is involved, from less than 1 kW to above 1 MW \cite{Rasmussen2005}. It is stated that about 30\% of the total energy all over the world is consumed by Heating, Ventilating, and Air Conditioning (HVAC) systems, as well as refrigerators and water heaters \cite{jahangeer2011numerical}. Therefore, the weight of the refrigeration processes on energy and economic balances is not in any way negligible \cite{buzelin2005experimental}. Indeed, regarding supermarkets and grocery stores, they are known to represent one of the largest consumers in the energy field, being 60\% of their consumption linked to refrigeration processes \cite{suzuki2011analysis}. Furthermore, it is reported that their average energy intensity is up to 500 kWh/m\textsuperscript{2} a year in USA, that corresponds to more than twice the energy consumed by a hotel or an office building \cite{US_EPA2}. 

Great effort has been made over the last decades to improve the overall energy efficiency of current refrigeration systems and reduce their environmental impact, through enhanced design of equipment (heat exchangers, compressors, valves, etc.), use of environmental-friendly refrigerants, and the application of advanced control and optimization strategies \cite{bejarano2017suboptimal,yin2018energy}. Furthermore, in recent years, a novel line of research regarding cold-energy management has been developed. Adding a thermal energy storage (TES) system to the canonical refrigeration cycle offers a number of advantages already exploited in thermal energy applications, for instance in distributed solar collector fields \cite{lima2016temperature,rubio2018optimal,navas2018optimal}. Since the TES system acts as an energy buffer, it is no longer necessary to produce exactly the required cooling demand at every moment. This feature allows to streamline the system capacity, in such a way that one can count both on the refrigeration cycle itself and on the cold energy already produced and stored in the TES system in order to address peak-demand periods. It also implies that the refrigeration cycle can work in more advantageous conditions to improve sustainability and efficiency. An additional advantage that arises from the decoupling of demand and production is the opportunity of scheduling the cooling generation to reduce the daily operating cost, considering the energy price and the market fluctuations (\emph{peak-shifting}) \cite{dincer2002bthermal,rismanchi2012energy}. This work is focused on the latter advantage, being the remaining ones more linked to the design stage.

Many commercial and under development solutions choose phase change materials (PCM) instead of sensible-heat ones for the TES system, as reported in some complete reviews \cite{mehling2008heat,oro2012review}. The main reason is related to thermodynamic properties, fitting better to energy storage in the case of PCM: higher heat capacity and minor temperature variations in latent state. In addition to the material, other differentiating factors between the diverse technologies are the encapsulation and the interface between the PCM and the heat transfer fluid (HTF), prevailing the packed bed technology over other layouts \cite{verma2008review,dutil2011review}.

Regarding cold-energy management, different strategies have been proposed in the literature. For instance, three works by Wang \emph{et al.}\cite{wang2007anovel,wang2007bnovel,wang2007cnovel} address the design, modelling, and control of a large-scale HVAC plant backed up by a ring of PCM-based TES tanks. Whereas the first work undertakes the system design and the second one addresses modelling, the third work by Wang \emph{et al.} presents a control strategy based on activation and deactivation of the TES tanks, where a measure of the overall performance is intended to be maximized. Moreover, Mossafa \emph{et al.}\cite{mosaffa2014advanced} rely on an exergy analysis to develop a management strategy consisting in combining alternatively the different PCM modules backing up the HVAC system. 

Other works apply techniques based on model predictive control (MPC) to energy management of TES-backed-up refrigeration systems in different applications. Shafiei \emph{et al.} \cite{shafiei2014model} propose a MPC strategy for a large-scale refrigeration plant, where the main objective is to track a given reference on the electric energy consumption. An optimization problem is posed, where the reference on the evaporating temperature is calculated as a virtual control variable, using an estimation on the energy stored and released from the TES tank. Moreover, Deng \emph{et al.} \cite{deng2015model} consider the optimal scheduling problem for a campus central plant equipped with a bank of multiple electrical chillers and a thermal energy storage. At each time step, the MPC problem is represented as a large-scale mixed-integer non-linear programming (MINLP) problem. In order to ensure the computational tractability of the problem, a suboptimal solution is proposed, where the optimal TES operation profile is obtained by solving a dynamic programming problem at every horizon, and the optimal chiller operations are obtained by solving a mixed integer linear programming (MILP) problem at every time step with a fixed TES operation profile.

MPC-based management strategies have also been applied to TES-backed-up refrigerated freight transport. For instance, Shafiei \emph{et al.} \cite{shafiei2015model} study a configuration where the TES tank is arranged in parallel with the refrigeration cycle; a prediction on the cooling demand is computed using the delivery route profile, traffic information, and weather forecast. Furthermore, another different configuration where the TES tank is arranged in series is analysed by Schalbart \emph{et al.} \cite{schalbart2015ice}. In this configuration, the refrigeration cycle charges the TES tank, that is in turn discharged while cooling down an ice-cream warehouse. A MPC-based strategy is applied to manage the system, where the prediction model includes a steady-state refrigeration cycle submodel and energy balances of the TES tank and the warehouse. The main objective is to ensure the product quality in long-term storage, assessed by means of a given ice-cream crystal size, while minimizing energy consumption.

\begin{figure}[h]
	\centerline{\includegraphics[width=7cm,trim = 40 110 35 155,clip]{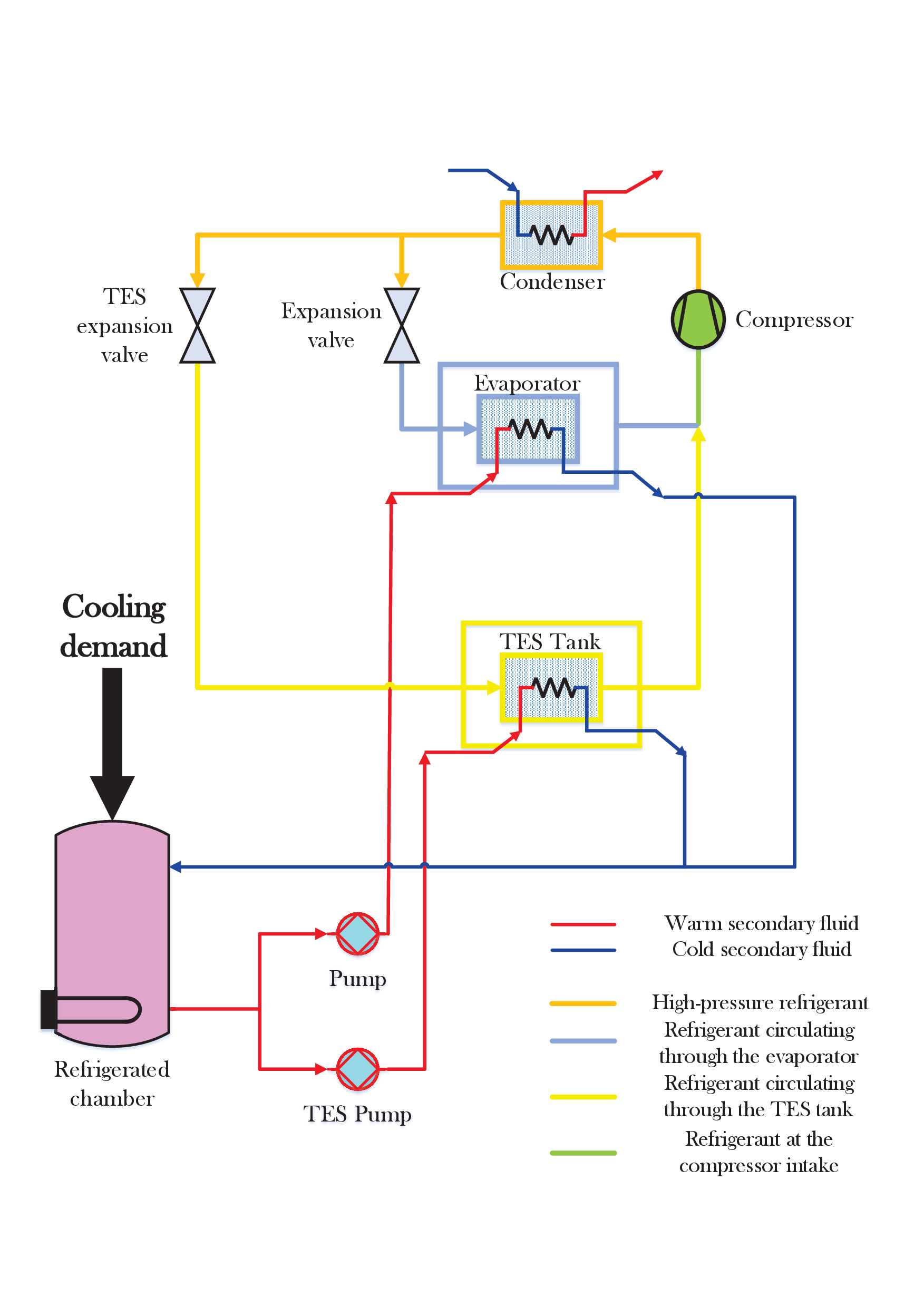}}
	\caption{Layout of the TES-backed-up refrigeration system considered in this work, where the TES tank is arranged in parallel with the evaporator.}
	\label{figEsquemaCicloPCM}
\end{figure}

The layout of the TES-backed-up refrigeration system considered in this work is represented in Figure \ref{figEsquemaCicloPCM}. In this case, the TES tank has been designed to be arranged in parallel with the evaporator of an existing refrigeration facility, whose original features have been described in previous works \cite{bejarano2015multivariable,GB_JE_2015}.

In the layout represented in Figure \ref{figEsquemaCicloPCM}, it is observed how the refrigerated chamber consists of a tank filled with a certain fluid (that will be called as secondary fluid from now on), that is cooled when circulating both through the evaporator and the TES tank, discharging the latter. An electric resistance is used in the refrigerated chamber to simulate the cooling demand. Concerning the refrigerant, it circulates through the canonical refrigeration cycle (compressor, condenser, expansion valve, and evaporator), but it may also circulate through the TES tank while charging it. This fact implies that the cold HTF (refrigerant) differs from the warm HTF (secondary fluid), unlike in the packed bed technology, where the same fluid is used as the cold HTF and the warm HTF. The TES layout was presented in a previous work by the authors \cite{bejarano2019NMPC} and it includes a number of PCM cylinders, two bundles of tubes that correspond to the refrigerant and the secondary fluid, and the so-called \emph{intermediate fluid} bathing all pipes and PCM cylinders. This setup is in turn very similar to that presented in the work by Bejarano \emph{et al.} \cite{Bejarano2017NovelSchemePCM}, being in this case the PCM encapsulation in the form of cylinders, instead of the spheres proposed in the aforementioned work. 

\begin{figure}[h]
	\centering
	\includegraphics[width=11.0cm,trim = 160 60 180 332,clip]
	{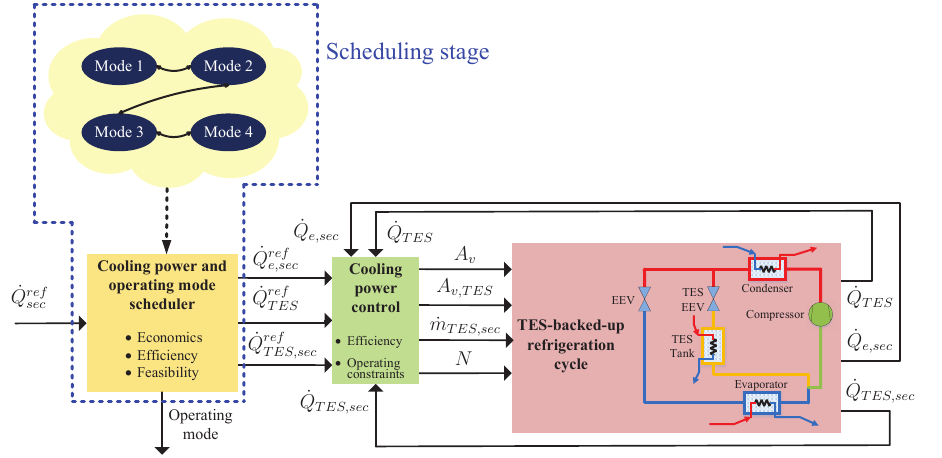}
	\caption{Scheduling and control strategy for the TES-backed-up refrigeration system \cite{bejarano2019NMPC}.}
	\label{figSchedulingControlStrategy}
\end{figure}

The addition of the TES tank to the canonical vapour-compression refrigeration cycle allows to decouple cooling demand and production, that might result in reducing the daily operating cost, as long as the cooling production is scheduled according to the energy market fluctuations. Furthermore, the combined refrigeration on the secondary fluid due to its circulation through the evaporator and the TES tank allows to satisfy peak demand that the original refrigeration cycle could not by itself. Obviously, the TES tank capacity must also be considered, in order to keep it in a latent state. Therefore, in order to satisfy a given demand profile that might require to provide some cooling power to the secondary fluid, both at the evaporator and the TES tank, as well as to reduce the overall operating cost, a challenging scheduling and control problem arises. The cascade strategy proposed in the work by Bejarano \emph{et al.} \cite{bejarano2019NMPC}, and shown in Figure \ref{figSchedulingControlStrategy}, is also considered in this work.

In the strategy described in Figure \ref{figSchedulingControlStrategy}, the scheduling stage is intended to compute the references on the main cooling powers involved in the problem: the one provided to the secondary fluid at the evaporator $\dot{Q}_{e,sec}^{ref}$, and the TES charging and discharging powers $\dot{Q}_{TES}^{ref}$ and $\dot{Q}_{TES,sec}^{ref}$, in such a way that the cooling demand $\dot{Q}_{sec}^{ref}$ is satisfied in real time and the operating cost is minimized. Then, a low-level cooling power control is applied to get the hybrid system to actually provide the required cooling powers, by driving the available control variables, namely the compressor speed $N$, the expansion valve openings $A_v$ and $A_{v,TES}$, and the TES pump through the virtual manipulated variable $\dot{m}_{TES,sec}$, that corresponds to the reference on the secondary mass flow circulating through the TES tank. Since the low-level cooling power control has been explained in detail in the aforementioned work by Bejarano \emph{et al.} \cite{bejarano2019NMPC}, this work is mostly focused on the scheduling stage, that is approached as a hybrid optimal control problem.

A scheduling strategy based on non-linear model predictive (NMPC) techniques has been proposed in the same aforementioned work by Bejarano \emph{et al.} \cite{bejarano2019NMPC}. In that strategy, the scheduling problem is posed as a non-linear optimization problem where the decision variables turn out to be the references on the three cooling powers involved throughout a given horizon. However, the operating mode scheduling is not included in the optimization problem, but it must be set off-line according to the predicted demand profile and the constraints on the achievable cooling powers. Moreover, some weights in the objective function must be manually tuned in such a way that, for instance, the TES tank is charged as much as necessary to satisfy a later peak demand. That is the main shortcoming of the NMPC-based scheduler, that is overcome in the scheduling strategy proposed in this work by including the operating mode scheduling within the optimization problem. Then, the decision set is not only comprised of the references on the relevant cooling powers, but some binary variables are also included, defining whether it is more suitable for the TES tank to be charged, discharged, or kept in a stand-by state all throughout the prediction horizon. This is the main contribution of the strategy proposed here, with respect to the NMPC-based scheduler described in the previous work by Bejarano \emph{et al.} \cite{bejarano2019NMPC}. 

Including discrete variables in the optimization problem, in addition to the continuous ones, introduces a new challenge from the point of view of the mathematical problem formulation and a corresponding suitable solving procedure. Several optimization algorithms have been proposed in the literature to control hybrid systems of variable configuration, where MPC turns out to be the most used strategy \cite{zengshan2001optimal,bemporad2002hybrid,hedlund2002convex,charbonnaud2003robust,mhaskar2005predictive,zambrano2010application}. MINLP solvers are frequently used to address this class of problems, that includes in this case some constraints ensuring the satisfaction of the cooling demand in real time, as well as those related to the limits on the achievable cooling powers. Furthermore, the TES \emph{charge ratio} $\gamma_{TES}$ is preferred to remain within a given security range that corresponds to the latent state. Therefore, a prediction model on how $\gamma_{TES}$ evolves when applying charging/discharging cooling power is needed. In the work by Bejarano \emph{et al.} \cite{bejarano2019NMPC} a simplified model focused on the dominant dynamics related to heat transfer within the TES tank was proposed, based on the highly time-efficient approach explained in a previous work by Bejarano \emph{et al.} \cite{bejarano2018efficient}. That model was used within the optimization as the prediction model. Nevertheless, it turns out to be too complex to be used within the proposed scheduler, since the MINLP solvers require simpler models to ensure reasonable solving time. 

In this paper an even simpler prediction model is used within the optimization, obtained by linearising the aforementioned simplified model. It is computed by using the main idea of a recent technique in predictive control, called Practical Non-linear Model Predictive Control (PNMPC) \cite{plucenio2007practical,Plucenio}. This technique seeks a linear representation of the predicted output with regard to the future control actions, but linearisation at a given equilibrium point is not considered. Instead, an approximation to the calculation of the forced response using the gradients, and therefore a first-order linearisation based on the corresponding Jacobian matrix, is recomputed at every sampling time.

State estimation is also addressed, given that the TES state vector is not fully measurable. In particular, the \emph{charge ratio} $\gamma_{TES}$ must be estimated according to measurable variables such as the temperature of the \emph{intermediate fluid}. Regarding the scheduling solving algorithm, a branch-and-bound method is applied \cite{lawler1966branch,christofides1987project,demeulemeester2000discrete}.

Some simulation results provided by the proposed scheduling strategy for a challenging cooling demand profile, that requires the combined use of the evaporator and the TES tank to face the peak demand, are presented and discussed. Moreover, the main advantages of the MINLP-based strategy over the aforementioned NMPC-based scheduler are remarked, concerning operating cost, ease of tuning, and ability to adapt to cooling demand variations. 

The brief is organised as follows. Section \ref{secSystemDescription_OpModes} provides some details about the TES tank and the operation of the hybrid system in the most useful operating modes. The proposed scheduling strategy is developed in Section \ref{secScheduling}, where the model linearisation based on the PNMPC ideas is detailed, as well as the constraints and the objective function imposed in the optimization. Section \ref{secEstimation} is devoted to the TES state estimator, inspired in the time-efficient discrete model developed by Bejarano \emph{et al.} \cite{bejarano2018efficient}. Section \ref{secCaseStudy} describes a case study for a challenging cooling demand profile, where some simulation results of the proposed strategy are discussed and its advantages over the NMPC-based scheduler are remarked. Finally, the main conclusions to be drawn and some future work are expressed in Section \ref{secConclusions}.

%%%%%%%%%%%%%%%%%%%%%%%%%%%%%%%%%%%%%%%%%%%%%%%%%%%%%%%%%%%%%%%%%%%%%%%%%%%%%%%%%%%%%

\section{System description and operating modes} \label{secSystemDescription_OpModes}

\subsection{Notation} \label{subSecNotation}

The notation followed throughout the work is detailed in Table \ref{tabSymbols}. 

\begin{table}[h]
	\centering
	\caption{Italic and greek symbols, as well as subscript/superscript notation}
	\label{tabSymbols}
	\begin{tabular} { L{1.2cm} L{6.0cm} L{1.8cm}  L{1.8cm} L{4.7cm} }
		\toprule
		\multicolumn{3}{c}{\textbf{Italic symbols}} & \multicolumn{2}{c}{\textbf{Subscripts}} \\
		\midrule
		\emph{Symbol} & \emph{Description} & \emph{Units} & \emph{Symbol} & \emph{Description} \\ 
		\midrule
		$A$ & Opening & \% & $coat$ & coating \\
		$c$ & Specific heat capacity & J kg\textsuperscript{-1} K\textsuperscript{-1} & $e$ & evaporator \\
		$D$ & Diameter & m & $forced$ & forced response \\
		$e$ & Thickness & m & $free$ & free response \\
		$\mathbold{G}$ & Dynamic matrix & -- & $in$ & inlet/input \\
		$h$ & Specific enthalpy & J kg\textsuperscript{-1} & $int$ & intermediate fluid \\
		$J$ & Objective function & \euro & $lay$ & cylindrical layer \\
		$k$ & Discrete step time & -- & $out$ & outlet/output \\
		$L$ & Length & m & $PNMPC$ & related to PNMPC \\
		$m$ & Mass & kg & $p$ & constant pressure \\
		
		$\dot m$ & Mass flow rate & g s\textsuperscript{-1} & $past$ & past and current \\
		$N$ & Compressor speed & Hz & $pcm$ & Phase Change Material \\
		$\mathbold{NLF}$ & Generic non-linear function & -- & $predict$ & predicted \\
		$n$ & Number of elements (e.g. PCM cylinders) & -- & $refr$ & refrigerant \\	
		$P$ & Pressure & Pa & $SH$ & superheating \\
		$PH$ & Prediction horizon & -- & $sec$ & secondary fluid \\
		$\dot Q$ & Cooling power & W & $surr$ & surroundings \\
		$T$ & Temperature & K & $tank$ & tank \\
		$t$ & Time & h & $TES$ & Thermal Energy Storage \\
		$U$ & Internal energy & J & $v$ & expansion valve \\
		$\mathbold{u}$ & Input vector & -- & & \\
		$V$ & Volume & m\textsuperscript{3} & & \\
		\cmidrule{4-5}
		$w$ & Weight in the objective function & \euro $\;$ W\textsuperscript{-1} & \multicolumn{2}{c}{\textbf{Superscripts}} \\
		\cmidrule{4-5}
		$\mathbold{x}$ & State vector & -- & \emph{Symbol} & \emph{Description} \\ 
		\cmidrule{4-5} 
		$\mathbold{y}$ & Output vector & -- & $lat$ & latent state \\
		\cmidrule{1-3}
		\multicolumn{3}{c}{\textbf{Greek symbols}} & $lat+$ & Maximum enthalpy latency point \\ 
		\cmidrule{1-3}  
		\emph{Symbol} & \emph{Description} & \emph{Units} & $lat-$ & Minimum enthalpy latency point \\
		\cmidrule{1-3}
		$\gamma$ & \emph{Charge ratio} & -- & $max$ & maximum \\
		$\delta$ & Binary variable & -- & $min$ & minimum \\
		$\kappa$ & Thermal conductivity & W m\textsuperscript{-1} K\textsuperscript{-1} & $ref$ & reference \\
		$\rho$ & Density & kg m\textsuperscript{-3} & & \\ 
		$\mathbold{\psi}$ & Partial decision set & -- & & \\
		$\mathbold{\Omega}$ & Decision set & -- & & \\
		\bottomrule
	\end{tabular}
\end{table}

\subsection{Hybrid system description} \label{subSecDescription}

As detailed in some previous works by the authors \cite{bejarano2019NMPC,Bejarano2017NovelSchemePCM,bejarano2018efficient}, an existing versatile two-compression-stage, two-load-demand experimental refrigeration plant located at the Department of Systems Engineering and Automatic Control at the University of Seville (Spain) is complemented by \emph{ad-hoc} designed TES tanks based on PCM. Although the plant can be configured to work with up to two compression stages and two evaporators, a reduced complexity configuration is first studied, including a canonical refrigeration cycle with only one compressor and one evaporator, where the TES tank is set up in parallel with the evaporator, as previously shown in Figure \ref{figEsquemaCicloPCM}. Additional elements have been deployed, such as the TES expansion valve and the TES pump, to drive the refrigerant and the secondary fluid, respectively, to the TES tank. Further information about the embedding of the TES tank in the original experimental facility can be found in the aforementioned related literature.

\begin{figure}[h]
	\centerline{\includegraphics[width=13cm,trim = 150 450 150 150,clip]
		{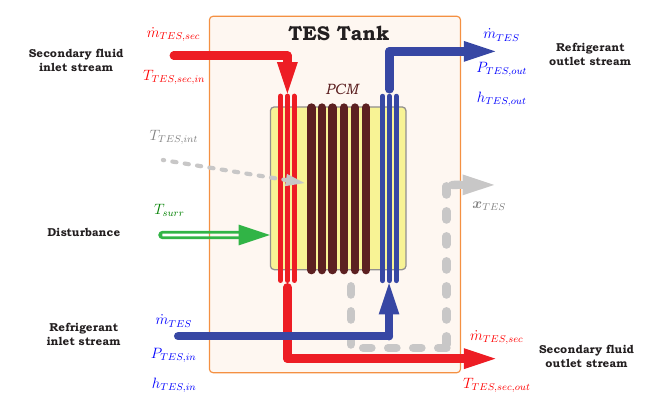}}
	\caption{TES tank setup.}
	\label{figTEStank_esquema}
\end{figure}

Figure \ref{figTEStank_esquema} shows the TES tank setup, where two bundles of pipes are deployed, corresponding to the warm HTF (secondary fluid) and the cold HTF (refrigerant). All pipes are bathed in the so-called \emph{intermediate fluid}, which presents high thermal conductivity and low heat capacity, while a counter-current configuration of the refrigerant and secondary fluid pipes is chosen to promote homogeneous heat transfer between the fluids running through the pipes and the \emph{intermediate fluid}. Furthermore, the PCM is confined in $n_{pcm}$ steel cylinders, also being dipped in the \emph{intermediate fluid}. The latter is completely still inside the TES tank, at constant atmospheric pressure, though closed. Its high thermal conductivity provides efficient heat transfer between the refrigerant and the PCM cylinders, during the charging cycle, or between the PCM cylinders and the secondary fluid, during the discharging cycle. Consequently, the \emph{intermediate fluid} is assumed to have, at each instant, a homogeneous temperature in the whole tank volume. Further information about the TES tank setup can be found in the related literature \cite{bejarano2019NMPC,Bejarano2017NovelSchemePCM,bejarano2018efficient}.

The main input and output variables of the TES tank can also be noticed in Figure \ref{figTEStank_esquema}, defining the mass flow and thermodynamic state description of the different inlet and outlet streams. In the case of the refrigerant, the \{$P$ -- $h$\} pair completely describes the thermodynamic state of both the inlet and outlet streams, whereas the temperature suffices in the case of the secondary fluid inlet/outlet streams. The ambient temperature $T_{surr}$ acts as a disturbance, while the temperature of the \emph{intermediate fluid} $T_{TES,int}$ is one of the TES tank state variables included in the state vector $\mathbold{x}_{TES}$,

\begin{equation}
	\mathbold{x}_{TES} = 
	\left[
		\begin{matrix} 
			h_{pcm,1} \\ 
			h_{pcm,2} \\ 
			\cdots \\ 
			h_{pcm,n_{lay}} \\ 
			T_{TES,int} 
		\end{matrix}
	\right] 
	\in \mathbb{R}^{n_{lay}+1} \; .
	\label{eq_x_TES}
\end{equation}

Moreover, as shown in \eqref{eq_x_TES}, $\mathbold{x}_{TES}$ includes the thermodynamic state of the $n_{lay}$ cylindrical layers into which every PCM cylinder is conceptually divided, represented by their specific enthalpy $h_{pcm,k} \; \forall k \in [1,n_{lay}]$, according to the discrete model proposed in the previous modelling works by Bejarano \emph{et al.} \cite{Bejarano2017NovelSchemePCM,bejarano2018efficient}.

The TES tank \emph{charge ratio} $\gamma_{TES}$ can be computed from the PCM enthalpy distribution as indicated in \eqref{eq_gamma_TES}, where it is assumed that each of the $n_{pcm}$ PCM cylinders within the TES tank presents, at any time, the same thermodynamic behaviour. Then, $U_{TES}^{max}$ ($U_{TES}^{min}$) corresponds to the maximum (minimum) latent thermal energy that can be stored in the whole TES tank in \eqref{eq_gamma_TES}, while $U_{pcm}^{max}$ ($U_{pcm}^{min}$) corresponds to the maximum (minimum) latent energy that can be stored in a single PCM cylinder. All these terms are constant, while $U_{TES} = n_{pcm}\,U_{pcm}$ refers to the variable latent energy stored in the whole TES tank. Then, $U_{pcm}$ can be computed from the enthalpy distribution within the PCM cylinders, considering the mass of every cylindrical layer ($m_{lay} = \frac{m_{pcm}}{n_{lay}}$) and its specific enthalpy $h_{pcm,k}\; \forall k \in [1,n_{lay}]$. Eventually, since it is interesting to measure the stored cold-thermal energy, instead of the mere thermal energy, the TES tank \emph{charge ratio} $\gamma_{TES}$ is defined as a normalised index between 0 and 1 in efficient storing conditions, as indicated below: 

\begin{equation}
	\begin{aligned}
		\gamma_{TES} &= \frac{U_{TES}^{max} - U_{TES}}{U_{TES}^{max} - U_{TES}^{min}} \;, \\
		U_{TES}^{max} &= n_{pcm} \, U_{pcm}^{max} \;, \\
		U_{TES}^{min} &= n_{pcm} \, U_{pcm}^{min} \;, \\
		U_{TES} &= n_{pcm} \, U_{pcm} \;, \\
		U_{pcm}^{max} &= V_{pcm}^{lat+} \, \rho_{pcm}^{lat+} \, h_{pcm}^{lat+} \;, \\
		U_{pcm}^{min} &= V_{pcm}^{lat-} \, \rho_{pcm}^{lat-} \, h_{pcm}^{lat-} \;, \\
		U_{pcm} &= \frac{m_{pcm}}{n_{lay}} \, \sum_{k=1}^{n_{lay}} h_{pcm,k} \;, \\
		\gamma_{TES} &= \frac{U_{pcm}^{max} - U_{pcm}}{U_{pcm}^{max} - U_{pcm}^{min}} \;. \\
	\end{aligned}
	\label{eq_gamma_TES}
\end{equation}
 
It is important to remark that $\gamma_{TES}$ is not a state variable, since a single value of $\gamma_{TES}$ can be obtained with different enthalpy distributions within the PCM cylinders.

\subsection{Operating modes} \label{subSecOperatingModes}

As stated in Section \ref{secIntroduction}, there are three main cooling powers generated in the system: 

\begin{itemize}
	\item The cooling power transferred from the refrigerant to the secondary fluid at the evaporator, denoted as $\dot{Q}_{e,sec}$.
	\item The cooling power transferred from the \emph{intermediate fluid} to the secondary fluid at the TES tank, denoted as $\dot{Q}_{TES,sec}$.
	\item The cooling power transferred from the refrigerant to the \emph{intermediate fluid} at the TES tank, denoted as $\dot{Q}_{TES}$.
\end{itemize}

\begin{figure}[h]
	\centering
	\includegraphics[width=13cm,trim = 140 510 160 180,clip]
	{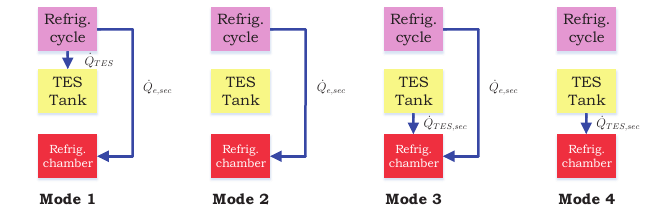}
	\caption{Most suitable operating modes regarding the scheduling problem \cite{bejarano2019NMPC}.}
	\label{figOperatingModes}
\end{figure}

The sum of the first two items represents the total cooling power provided to the secondary fluid, that must match the cooling demand at any time, corresponding the second one to the TES discharging power. The last one represents the TES charging power. Up to eight operating modes can be defined according to all possible combinations of these three cooling powers; they all have been described and discussed in a previous work by Bejarano \emph{et al.} \cite{bejarano2019NMPC}. However, some combinations might be meaningless or not very useful for the problem of satisfying a realistic cooling demand. The most suitable operating modes regarding the scheduling problem are modes 1 to 4, graphically described in Figure \ref{figOperatingModes}. 

Since it is assumed that a non-zero cooling demand must be satisfied all throughout the day, the refrigerated chamber must take at least one cooling power contribution, either provided by the refrigeration cycle at the evaporator, or supplied by the TES tank, or both. As long as the demand is attainable only by providing only cooling power at the evaporator, modes 1 and 2 can be scheduled: the only difference between them lies in whether the TES tank is simultaneously being charged or not, which is to be decided according to the \emph{charge ratio} and the demand forecast. If the latter is high enough to require the double cooling power contribution at the evaporator and at the TES tank, mode 3 must be scheduled. Eventually, mode 4 might be used when the demand is satisfiable just by discharging the TES tank, and this operation is economically advantageous, provided that the \emph{charge ratio} is high enough to satisfy the demand during the proposed period.

%%%%%%%%%%%%%%%%%%%%%%%%%%%%%%%%%%%%%%%%%%%%%%%%%%%%%%%%%%%%%%%%%%%%%%%%%%%%%%%%%%%%%

\section{Scheduling strategy} \label{secScheduling}

\subsection{Overview} \label{subSecOverview}

\begin{figure}[h]
	\centering
	\includegraphics[width=14.0cm,trim = 230 160 230 215,clip]
	{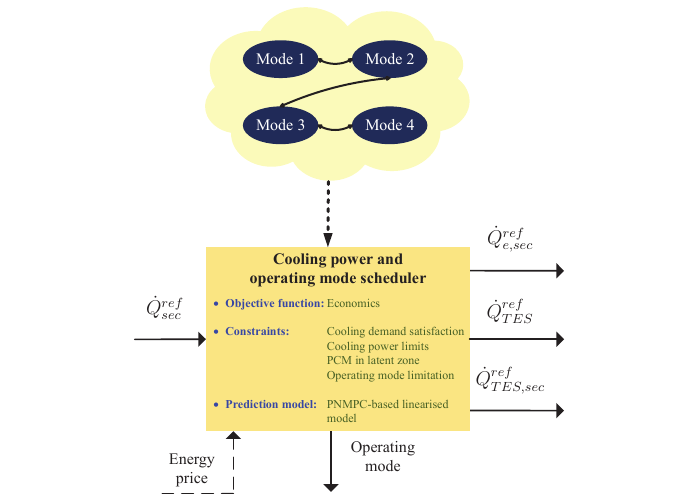}
	\caption{Proposed scheduling strategy for the TES-backed-up refrigeration system.}
	\label{figSchedulingStrategy}
\end{figure}

As stated in Section \ref{secIntroduction}, this work is focused on the scheduling stage of the cascade control strategy shown in Figure \ref{figSchedulingControlStrategy}, that is detailed in Figure \ref{figSchedulingStrategy}. The scheduling problem is posed as a mixed integer non-linear optimization problem, where a receding horizon strategy is applied. Therefore, given a certain prediction/control horizon, the objective of the scheduler is to compute some feasible references for $\dot{Q}_{e,sec}^{ref}$, $\dot{Q}_{TES}^{ref}$, and $\dot{Q}_{TES,sec}^{ref}$, in such a way that the cooling demand $\dot{Q}_{sec}^{ref}$ is satisfied at any time, the \emph{charge ratio} $\gamma_{TES}$ remains within a given range corresponding to the PCM latent zone, and the operating cost is minimized, according to the variable energy price that corresponds to actual market fluctuations. The mixed features of the optimization arise from the inclusion of binary variables that define the operating mode scheduling and are part of the decision variable set. 

In the following subsections the decision variable set, the constraints, the prediction model, and the objective function considered in the optimization problem are detailed. 

\subsection{Decision set} \label{subSecDecisionSet}

Firstly, the decision variable set is described. It comprises the references on two of the three relevant cooling powers throughout the prediction/control horizon $PH$: \{$\dot{Q}_{TES}^{ref}(t-1+k)$, $\dot{Q}_{TES,sec}^{ref}(t-1+k)$\} $\; \forall k \in [1,PH]$. As shown later, the references on the cooling power provided at the evaporator, $\dot{Q}_{e,sec}^{ref}(t-1+k) \; \forall k \in [1,PH]$, are computed from the decision variables by imposing the cooling demand satisfaction constraint. Moreover, two binary variables \{$\delta_{TES}(t-1+k)$, $\delta_{TES,sec}(t-1+k)$\} $\; \forall k \in [1,PH]$ are also included in the decision variable set, indicating whether the corresponding cooling powers are active or not all throughout the prediction horizon. Then, the decision variable set $\mathbold{\Omega} (t)$ is the one described by:

\begin{equation}
	\begin{aligned}
		\mathbold{\psi}(t-1+k) &= 
		\left[
			\begin{matrix}
				\dot{Q}_{TES}^{ref}(t-1+k) \\
				\dot{Q}_{TES,sec}^{ref}(t-1+k) \\
				\delta_{TES}(t-1+k) \\
				\delta_{TES,sec}(t-1+k) \\
			\end{matrix}
		\right] \qquad \forall k \in [1,PH] \;, \\
		\mathbold{\Omega} (t) &= 
		\left[
			\begin{matrix}
				\mathbold{\psi}(t) \\
				\mathbold{\psi}(t+1) \\
				\cdots \\
				\mathbold{\psi}(t-1+PH)  \\
			\end{matrix}
		\right] \;,
	\end{aligned}
	\label{eq_Decision_Set}
\end{equation}

\noindent where the partial decision set $\mathbold{\psi}$ includes the decision variables corresponding to a given instant $t-1+k \;\; \forall k \in [1,PH]$. As shown in \eqref{eq_Decision_Set}, the total number of decision variables in $\mathbold{\Omega} (t)$ is $4 \cdot PH$, among which half are continuous and half binary. It is important to note that a single horizon $PH$ is considered, in such a way that the control horizon matches the prediction one.

\subsection{Constraints} \label{subSecConstraints}

As shown in Figure \ref{figSchedulingStrategy}, the constraints imposed on the decision variables are related to several issues. First of all, the cooling demand must be satisfied at any time throughout the horizon. It is achieved by computing the references on the cooling power provided at the evaporator, $\dot{Q}_{e,sec}^{ref}(t-1+k) \; \forall k \in [1,PH]$, as indicated below:

\begin{equation}
	\dot{Q}_{e,sec}^{ref}(t-1+k) = \dot{Q}_{sec}^{ref}(t-1+k) - \dot{Q}_{TES,sec}^{ref}(t-1+k) \geq 0 \quad \forall k \in [1,PH] \;.
	\label{eq_Cooling_Demand_Constraint}
\end{equation}

Moreover, an auxiliary binary variable is defined, $\delta_{e,sec}(t-1+k) \; \forall k \in [1,PH]$, which indicates whether the corresponding cooling power $\dot{Q}_{e,sec}^{ref}(t-1+k)$ is active or not. It is computed according to the value resulting from the application of \eqref{eq_Cooling_Demand_Constraint} all throughout the horizon $PH$. Therefore, this formulation of the cooling demand satisfaction constraint allows to reduce the optimization problem size, given that the references on the cooling power provided at the evaporator and the corresponding binary variables are no longer included in the decision set $\mathbold{\Omega}$.

Secondly, the references on the cooling powers must be feasible and achievable by the TES-backed refrigeration systems when operating in all modes. The cooling power feasibility is imposed through the constraints:

\begin{equation}
		\left.
		\begin{aligned}
			\begin{matrix}
				\dot{Q}_{e,sec}^{ref}(t-1+k) \geq \delta_{e,sec}(t-1+k) \, \dot{Q}_{e,sec}^{ref,min}(t-1+k) \\
				\dot{Q}_{e,sec}^{ref}(t-1+k) \leq \delta_{e,sec}(t-1+k) \, \dot{Q}_{e,sec}^{ref,max}(t-1+k) \\
			\end{matrix} \\
			\begin{matrix}
				\dot{Q}_{TES}^{ref}(t-1+k) \geq \delta_{TES}(t-1+k) \, \dot{Q}_{TES}^{ref,min}(t-1+k) \\
				\dot{Q}_{TES}^{ref}(t-1+k) \leq \delta_{TES}(t-1+k) \, \dot{Q}_{TES}^{ref,max}(t-1+k) \\
			\end{matrix} \\
			\begin{matrix}
				\dot{Q}_{TES,sec}^{ref}(t-1+k) \geq \delta_{TES,sec}(t-1+k) \, \dot{Q}_{TES,sec}^{ref,min}(t-1+k) \\
				\dot{Q}_{TES,sec}^{ref}(t-1+k) \leq \delta_{TES,sec}(t-1+k) \, \dot{Q}_{TES,sec}^{ref,max}(t-1+k) \\
			\end{matrix} \\
		\end{aligned}
	\right\} \quad \forall k \in [1,PH]	\;.
	\label{eq_Feasibility_Constraints}
\end{equation} 

The power limit formulation shown in \eqref{eq_Feasibility_Constraints} allows to impose that a certain cooling power must be zero if the corresponding binary variable has been set to zero/false. Furthermore, the maximum and minimum values of the feasible cooling powers at every $k$ may depend not only on the operating mode defined by the binary variable set \{ $\delta_{e,sec}(t-1+k)$, $\delta_{TES}(t-1+k)$, $\delta_{TES,sec}(t-1+k)$ \}, but also on the \emph{charge ratio} $\gamma_{TES}$ and, what is more, on the specific enthalpy distribution inside the PCM cylinders, given by the state vector $\mathbold{x}_{TES}$ shown in \eqref{eq_x_TES}. This issue has been discussed in the work by Bejarano \emph{et al.} \cite{bejarano2019NMPC}: the main reason lies in the fact that, as the PCM cylinders are charged/discharged, the thermal resistance caused by the cylindrical shell in the sensible zone becomes greater, modifying the minimum and maximum achievable cooling powers. Moreover, since a typical TES tank operation is expected to schedule several partial charging/discharging processes, multiple moving freezing/melting boundaries are very likely to be present at the same time inside the PCM cylinders \cite{Bejarano2017NovelSchemePCM,bejarano2018efficient}. Indeed, it is only the outermost cylindrical shell in the sensible zone that defines the applicable thermal resistance, and therefore it depends on the PCM cylinder \emph{history}, and that is why an estimation of the enthalpy distribution inside the PCM cylinders is needed to impose the right limits on the achievable cooling powers.   

Thirdly, it is intended that the PCM remains in the latent zone throughout the complete horizon, namely $\gamma_{TES}$ must remain within the range [0, 1]. However, security limits $\gamma_{TES}^{min}>0$ and $\gamma_{TES}^{max}<1$ are usually imposed as indicated below:

\begin{equation}
	\begin{aligned}
		\left.
			\begin{matrix}
				\gamma_{TES}(t-1+k) &\geq \gamma_{TES}^{min} \\
				\gamma_{TES}(t-1+k) &\leq \gamma_{TES}^{max} \\
			\end{matrix}
		\right\} \quad \forall k \in [1,PH] \;.
	\end{aligned}
	\label{eq_GammaTESConstraints}
\end{equation}

Some constraints may also be additionally imposed on the binary variables, in order to ensure that the scheduled operating modes correspond to those analysed in subsection \ref{subSecOperatingModes}. These constraints, expressed as logical conditions, 

\begin{equation}
	\begin{aligned}
		\left.
			\begin{matrix}
				\delta_{e,sec}(t-1+k) \; || \; \delta_{TES,sec}(t-1+k) &= \textnormal{true} \\
				\delta_{TES}(t-1+k) \; \& \; \delta_{TES,sec}(t-1+k) &= \textnormal{false} \\
			\end{matrix}
		\right\} \quad \forall k \in [1,PH] \;,
	\end{aligned}
	\label{eq_BinaryConstraints}
\end{equation}

\noindent imply that only operating modes 1 to 4 can be scheduled. Indeed, the first condition shown in \eqref{eq_BinaryConstraints} forces that there always exists at least one contribution to the satisfaction of the cooling demand (at the evaporator and/or at the TES tank), while the second condition ensures that no simultaneous TES charging and discharging is allowed.

\subsection{Prediction model} \label{subSecPredictionModel}

The constraints imposed on $\gamma_{TES}$ involve considering a prediction model of the TES tank behaviour in the optimization problem. Given that the simplified dynamic model proposed in the work by Bejarano \emph{et al.} \cite{bejarano2019NMPC} turns out to be too complex to be included in the mixed-integer optimization, a first-order linearised model calculated at every sampling time without using the equilibrium point concept is used as the prediction model. This idea has been proposed within the predictive control technique called PNMPC \cite{plucenio2007practical,Plucenio} to obtain a linear representation of the predicted output with regard to the future control actions. 

The simplified, non-linear model, focused on the dominant dynamics related to heat transfer within the TES tank, in state-space form, is shown below:

\begin{subequations}
	\begin{equation}
		\begin{aligned}
			\mathbold{x}_{TES} &= 
			\left[
				\begin{matrix} 
					h_{pcm,1} \\ 
					h_{pcm,2} \\ 
					\cdots \\ 
					h_{pcm,n_{lay}} \\ 
					T_{TES,int} \\
				\end{matrix}
			\right] 
			\in \mathbb{R}^{n_{lay}+1} \;, \\
			\mathbold{u}_{TES} &= 
			\left[
				\begin{matrix} 
					\dot{Q}_{TES}^{ref} \\ 
					\dot{Q}_{TES,sec}^{ref} \\ 
				\end{matrix}
			\right] 
			\in \mathbb{R}^{2} \;, \\
			\mathbold{y}_{TES} &= 
			\left[
				\begin{matrix} 
					\Delta \gamma_{TES} \\ 
					\Delta T_{TES,int} \\ 
				\end{matrix}
			\right] 
			\in \mathbb{R}^{2} \;, \\
		\end{aligned} 
	\end{equation}
	\begin{equation}
		\begin{aligned}
			\left.
				\begin{matrix}
					\begin{aligned}
						\mathbold{x}_{TES} (t+k) &= \mathbold{f}(\mathbold{x}_{TES} (t-1+k),\mathbold{u}_{TES} (t+k-1)) \\
						\mathbold{y}_{TES} (t+k-1) &= \mathbold{g}(\mathbold{x}_{TES} (t-1+k)) \\
					\end{aligned}
				\end{matrix}
			\right\}
			\quad \forall k \in [1,PH] \; ,
		\end{aligned} 
	\end{equation}
	\label{eq_Simplified_Model}			
\end{subequations}

\noindent where $\mathbold{x}_{TES}$ refers to the TES tank state vector, $\mathbold{u}_{TES}$ is the TES tank input vector, and $\mathbold{y}_{TES}$ corresponds to the output vector, that may include the increments on the temperature of the \emph{intermediate fluid} $\Delta T_{TES,int}$ and the \emph{charge ratio} $\Delta \gamma_{TES}$, defined as:

\begin{equation}
	\begin{aligned}
		\left.
			\begin{matrix}
				\begin{aligned}
					\Delta \gamma_{TES}(t+k) &= \gamma_{TES}(t+k) - \gamma_{TES}(t-1+k) \\
					\Delta T_{TES,int}(t+k) &= T_{TES,int}(t+k) - T_{TES,int}(t-1+k) \\
				\end{aligned}
			\end{matrix}
		\right\}
		\quad \forall k \in [1,PH] \; .
	\end{aligned} 
	\label{eq_IncrementDefinitions}		
\end{equation}

The main idea behind the PNMPC formulation is to deal with non-linear systems using the MPC techniques developed for linear systems \cite{plucenio2007practical,Plucenio}. In conventional linear MPC techniques, the vector of predicted outputs $\mathbold{y}_{predict}$ can be expressed as a linear function of the vector of future control inputs $\mathbold{u}$, where the free response $\mathbold{y}_{free}$ and the forced response $\mathbold{y}_{forced}$ are explicitly separated, being $\mathbold{G}$ a constant matrix denominated dynamic matrix of the model, as shown below:

\begin{equation}
	\begin{aligned}
		\mathbold{y}_{predict} = \mathbold{y}_{free} + \mathbold{y}_{forced} = \mathbold{y}_{free} + \mathbold{G} \, \mathbold{u} \; . \\
	\end{aligned}
	\label{eq_DMC_GPC}
\end{equation}  

The system shown in \eqref{eq_Simplified_Model} can also be expressed as

\begin{equation}
	\begin{aligned}
		\mathbold{y}_{predict} = \mathbold{NLF}(\mathbold{y}_{past},\mathbold{u}_{past},\mathbold{u}) \;, \\
	\end{aligned}
	\label{eq_nonlinear}
\end{equation}  

\noindent where the predicted output vector $\mathbold{y}_{predict}$ turns out to be a certain non-linear function of the current and past outputs $\mathbold{y}_{past}$, the past control inputs $\mathbold{u}_{past}$, and the future control actions $\mathbold{u}$, where $\mathbold{NLF}$ refers to the arbitrary non-linear function that defines the system.

Following the structure of the linear MPC shown in \eqref{eq_DMC_GPC}, the predicted output vector $\mathbold{y}_{predict}$ can be divided in two parts: the free response $\mathbold{y}_{free}$ (only due to the current and past outputs $\mathbold{y}_{past}$ and the past control inputs $\mathbold{u}_{past}$), and the forced response $\mathbold{y}_{forced}$, affected by $\mathbold{u}$. Regarding $\mathbold{y}_{free}$, this variable is computed by applying zero future control actions to the original non-linear model, as indicated below:

\begin{equation}
	\begin{aligned}
		\mathbold{y}_{free} &= \mathbold{NLF}(\mathbold{y}_{past},\mathbold{u}_{past},\mathbold{u}=\bm 0) \;, \\
		\mathbold{y}_{forced} &\approx \mathbold{G}_{PNMPC} \,  \mathbold{u} \;, \\
		\mathbold{G}_{PNMPC} &= \dfrac{\partial \mathbold{y}_{predict}}{\partial \mathbold{u}} \biggr\rvert_{\mathbold{u} = \bm 0} \; .
	\end{aligned}
	\label{eq_GPNMPC}
\end{equation}  

Concerning $\mathbold{y}_{forced}$, an approximation consisting of a first-order linearisation of the MacLaurin series is proposed, since it is computed around $\mathbold{u} = \bm 0$, as described in \eqref{eq_GPNMPC}.

The matrix $\mathbold{G}_{PNMPC}$ represents the Jacobian matrix, namely the gradient of $\mathbold{y}_{predict}$ with respect to future control inputs. The numerical algorithm presented in the related literature for multiple-input-multiple-output (MIMO) systems is used to compute $\mathbold{G}_{PNMPC}$ \cite{plucenio2007practical,Plucenio,bejarano2017optimization}. State feedback is used at every sampling time when applying the non-linear model shown in \eqref{eq_Simplified_Model} and \eqref{eq_nonlinear} to avoid offset and close the loop. 

Once the Jacobian matrix $\mathbold{G}_{PNMPC}$ is calculated, together with the free response $\mathbold{y}_{free}$, a linear prediction model is available. This Jacobian matrix is used within the optimization procedure to obtain the predicted values of the output vector $\mathbold{y}_{TES}$: the predicted \emph{charge ratio} $\hat{\gamma}_{TES}$, whose predictions are mandatory to impose the constraints described in \eqref{eq_GammaTESConstraints}, and the predicted temperature of the \emph{intermediate fluid} $\hat{T}_{TES,int}$.  

However, it has been stated in subsection \ref{subSecConstraints} that, regarding the feasibility constraints shown in \eqref{eq_Feasibility_Constraints}, the maximum and minimum values of the achievable cooling powers depend not only on the \emph{charge ratio} $\gamma_{TES}$, but also on the enthalpy distribution inside the PCM cylinder, described by the TES tank state vector $\mathbold{x}_{TES}$. This fact implies that at least an estimate of where the outermost cylindrical layer in latent zone is located is needed to determine the applicable limits on the cooling powers, according to the operating mode defined by the binary decision variables. The PNMPC-based strategy used to obtain the prediction model is not suitable, but the estimated $\hat{\gamma}_{TES}$ can be used to compute an estimation of the cold energy transferred by every PCM cylinder during every sampling time within the prediction horizon $\Delta\hat{U}_{pcm}(t+k) \quad \forall k \in [1,PH]$, as indicated below:

\begin{equation}
	\begin{aligned}
		\Delta\hat{U}_{pcm}(t+k) = \Delta \gamma_{TES}(t+k) \, (U_{pcm}^{max} - U_{pcm}^{min}) = \Big[ \hat{\gamma}_{TES}(t+k) - \hat{\gamma}_{TES}(t-1+k) \Big] \, (U_{pcm}^{max} - U_{pcm}^{min})  \quad \quad \forall k \in [1,PH] \;. \\ 
	\end{aligned}
	\label{eq_U_TES}
\end{equation}  

Once estimated the cold energy transferred during every sampling time $\Delta\hat{U}_{pcm}$, and given the TES tank state at the initial point of the prediction horizon $\mathbold{\hat{x}}_{TES}(t) \equiv \mathbold{x}_{TES}(t)$, a recursive algorithm is applied to compute an estimation on the enthalpy distribution inside the PCM cylinder all throughout the prediction horizon. This algorithm is inspired in that proposed in the modelling work by Bejarano \emph{et al.} \cite{bejarano2018efficient} implementing the simplified dynamic model shown in \eqref{eq_Simplified_Model}. For every sampling time within the prediction horizon, $\forall k \in [1,PH]$, the algorithm is expressed as a step-by-step sketch, as follows:

\begin{enumerate}
	
	\item
	Starting from a given estimated state of the layered PCM cylinder and the \emph{intermediate fluid} $\mathbold{\hat{x}}_{TES}(t-1+k)$, an inward scanning sequence is performed, looking for the outermost layer $j_0$ in latent zone:
	
	\begin{equation}
		j_0 = \operatorname{max}\left\{j \in \{1,\ldots n_{lay}\} \mid
		h_{pcm}^{lat-} < \hat{h}_{pcm,j}(t-1+k) < h_{pcm}^{lat+} \right\} \;.
		\label{eq_j0}
	\end{equation}
	
	\item 
	Given the cold energy transferred $\Delta\hat{U}_{pcm}(t+k)$, the specific enthalpy of layer $j_0$ is updated accordingly:
	
	\begin{equation}
		\hat{h}_{pcm,j_0}(t+k) = \hat{h}_{pcm,j_0}(t-1+k) + \frac{\Delta\hat{U}_{pcm}(t+k) }{\hat{\rho}_{pcm,j_0}(t-1+k)\,\hat{V}_{pcm,j_0}(t-1+k)} \;,
		\label{eq_h_pcm_j0_forward}
	\end{equation}
	
	where $\hat{\rho}_{pcm,j_0}(t-1+k)$ and $\hat{V}_{pcm,j_0}(t-1+k)$ refer to the density and volume of layer $j_0$, computed from the estimated enthalpy $\hat{h}_{pcm,j_0}(t-1+k)$.
	
	\item 
	At this point, two possibilities arise:
	
	\begin{enumerate}[a)]
		
		\item
		Layer $j_0$ remains in the latent zone: $h_{pcm}^{lat-} < \hat{h}_{pcm,j_0}(t+k) < h_{pcm}^{lat+}$. That means that there is no change in the enthalpic state of the layers interior to $j_0$:
		
		\begin{equation}
			\hat{h}_{pcm,j}(t+k) = h_{pcm,j}(t-1+k) \qquad \forall j<j_0 \;.
			\label{eq_h_pcm_j_forward}
		\end{equation}
		
		Furthermore, layers exterior to $j_0$ are in the sensible zone. Depending on the sign of $\Delta\hat{U}_{pcm}(t+k)$, the estimation on their enthalpy is saturated to $h_{pcm}^{lat-}$ or $h_{pcm}^{lat+}$:
		
		\begin{equation}
			\begin{aligned}
				\hat{h}_{pcm,j}(t+k) = h_{pcm}^{lat-} \qquad \forall\, j>j_0 \quad \textnormal{if} \; \Delta\hat{U}_{pcm}(t+k) > 0  \;, \\
				\hat{h}_{pcm,j}(t+k) = h_{pcm}^{lat+} \qquad \forall\, j>j_0 \quad \textnormal{if} \; \Delta\hat{U}_{pcm}(t+k) < 0  \;. \\
			\end{aligned}	
			\label{eq_h_pcm_sensible}
		\end{equation}
		
		Then, the estimated state vector $\mathbold{\hat{x}}_{TES}(t+k)$ is computed as shown below:
		
		\begin{equation}
			\mathbold{\hat{x}}_{TES}(t+k) = 
			\left[
				\begin{matrix} 
					\hat{h}_{pcm,1}(t+k) \\ 
					\hat{h}_{pcm,2}(t+k) \\ 
					\cdots \\ 
					\hat{h}_{pcm,n_{lay}}(t+k) \\ 
					\hat{T}_{TES,int} (t+k)
				\end{matrix}
			\right] \;, 
			\label{eq_x_TES_estimated}
		\end{equation}
		
		 where $\hat{T}_{TES,int} (t+k)$ is obtained from \eqref{eq_IncrementDefinitions}.

		\item Layer $j_0$ quits the latent zone. That means that the latent energy of layer $j_0$ depleted some time before the sampling time expired, $\Delta t_{j_0} \leq \Delta t$, when the layer entered sensible zone. To continue with the algorithm, the cold energy transferred $\Delta\hat{U}_{pcm}(t+k)$ is updated as shown below, depending on the sign of $\Delta\hat{U}_{pcm}(t+k)$, given that a part of the original energy $\Delta\hat{U}_{pcm}(t+k)$ has been already transferred to layer $j_0$:
		
		\begin{equation}
			\begin{aligned}
				\Delta\hat{U}_{pcm}(t+k) = \Delta\hat{U}_{pcm}(t+k) - \left[\hat{h}_{pcm,j}(t-1+k) - h_{pcm}^{lat-}\right]\,\hat{\rho}_{pcm,j_0}(t-1+k)\,\hat{V}_{pcm,j_0}(t-1+k) \\
				\textnormal{if} \; \Delta\hat{U}_{pcm}(t+k) > 0  \;, \\ 
				\\
				\Delta\hat{U}_{pcm}(t+k) = \Delta\hat{U}_{pcm}(t+k) + \left[h_{pcm}^{lat+} - \hat{h}_{pcm,j}(t-1+k)\right]\,\hat{\rho}_{pcm,j_0}(t-1+k)\,\hat{V}_{pcm,j_0}(t-1+k) \\ 
				\textnormal{if} \; \Delta\hat{U}_{pcm}(t+k) < 0  \;. \\
			\label{eq_DeltaU_pcm_update}
			\end{aligned}	
		\end{equation}
		
		Then, the next inner layer, $j_0\!-\!1$, is established as the new outermost layer in latent zone, and the sequence restarts from step 2, applying the updated value of $\Delta\hat{U}_{pcm}(t+k)$ computed in \eqref{eq_DeltaU_pcm_update}.
		
	\end{enumerate}
	
\end{enumerate}

This algorithm allows to have a cold-energy-based estimation on where the outermost cylindrical layer in latent zone is located, obtained from the information provided by the linearised prediction model. It is important to note that no hypothesis about heat transfer between the PCM cylinder and the \emph{intermediate fluid} is considered, in such a way that the estimator is only based on the predicted \emph{charge ratio}, as well as on some thermodynamic properties of the PCM and geometric features of the TES tank. The proposed estimator is expected to be accurate enough, given that these properties and features are usually accurately known. The location of the outermost cylindrical layer in latent zone allows to impose the right cooling power limits on the feasibility constraints indicated in \eqref{eq_Feasibility_Constraints}, taking into account the thermal resistance caused by the cylindrical shell in the sensible zone, both during charging and discharging processes.

As the TES tank is charged/discharged, the cylindrical shell in the sensible zone grows and the minimum/maximum achievable charging/discharging cooling power is reduced. In order to ensure that the limits imposed on the feasibility constraints shown in \eqref{eq_Feasibility_Constraints} are actually achievable by the cycle during the complete optimization sampling time, the minimum value imposed is that corresponding to the position of the outermost cylindrical layer in the latent state at instant $t+k-1 \; \forall k \in [1,PH]$, whereas the maximum value corresponds to the position of the outermost cylindrical layer in the latent state at instant $t+k \; \forall k \in [1,PH]$. Nevertheless, when a transition between charging/discharging processes happens, the position of the outermost layer in the latent state is reset to the PCM cylinder edge, ans thus the minimum value is that corresponding to this situation, while the maximum value is computed according to the position of the outermost cylindrical layer in the latent state at instant $t+k \; \forall k \in [1,PH]$. The transition between charging/discharging processes is detected by comparing the predicted value of the temperature of the \emph{intermediate fluid} $\hat{T}_{TES,int}$, given by the linearised prediction model, with the phase-change temperature $T_{pcm}^{lat}$.

\subsection{Objective function} \label{subSecObjectiveFunction}

The objective function $J$, expressed as

\begin{equation}
	\begin{split}
		J = \sum_{k=1}^{PH} \; &w_{e,sec}(t-1+k) \; \dot{Q}_{e,sec}^{ref}(t-1+k) \; + \\ 
		+ \; &w_{TES}(t-1+k) \; \dot{Q}_{TES}^{ref}(t-1+k) \; + \\ 
		+ \; &w_{TES,sec}(t-1+k) \; \dot{Q}_{TES,sec}^{ref}(t-1+k) \;,
	\end{split}
	\label{eq_Objective_Function}
\end{equation}

\noindent includes only terms related to economic cost of cooling power generation all throughout the horizon $PH$. Note that in \eqref{eq_Objective_Function} the weights in the objective function $J$ of the cooling powers $\dot{Q}_{e,sec}^{ref}(t-1+k)$ and $\dot{Q}_{TES}^{ref}(t-1+k) \; \forall k \in [1,PH]$ correspond to the economic cost of producing such powers by the enhanced refrigeration cycle. However, the weights $w_{TES,sec}(t-1+k)$ are in this case set to zero since the fact of discharging the cold energy previously stored in the TES tank does not involve instantaneous cooling power production.

%%%%%%%%%%%%%%%%%%%%%%%%%%%%%%%%%%%%%%%%%%%%%%%%%%%%%%%%%%%%%%%%%%%%%%%%%%%%%%%%%%%%%

\section{State estimation} \label{secEstimation}

It has been stated in Section \ref{secScheduling} that the state estimation, $\hat{x}_{TES}$, is required, not only for feedback purposes, but also for enabling the estimation algorithm on where the outermost cylindrical layer in the latent state is located, in order to impose the right cooling power limits on the feasibility constraints described in \eqref{eq_Feasibility_Constraints}. Actually, the state vector $\mathbold{x}_{TES}$ is not fully measurable, since it includes specific enthalpies of the different PCM cylinder layers, as shown in \eqref{eq_Simplified_Model}. Although the temperature of the \emph{intermediate fluid} $T_{TES,int}$ is measurable, it will be used as the system feedback to develop the state estimator described in this section.

As stated in the modelling work by Bejarano \emph{et al.} \cite{Bejarano2017NovelSchemePCM,bejarano2018efficient}, the energy balance on the \emph{intermediate fluid} is given by 

\begin{equation}
	\begin{aligned}
		c_{p,TES,int} \, m_{TES,int} \, \Delta T_{TES,int} = \int_{0}^{\Delta t} \Big(\dot{Q}_{TES,sec} - \dot{Q}_{TES} - \dot{Q}_{TES,int} + \dot{Q}_{surr} \Big) dt \;, \\
	\end{aligned}
	\label{eq_EnergyBalance}
\end{equation}

\noindent where $\dot{Q}_{TES,int}$ refers to the cooling power transferred from the PCM cylinders to the \emph{intermediate fluid} (positive during discharging processes, negative when charging), $\dot{Q}_{surr}$ corresponds to thermal losses, and $\Delta t$ represents the sampling time.

According to the work by Bejarano \emph{et al.} \cite{bejarano2019NMPC}, the separation between the time scales of the scheduler and cooling power controller allows to assume that the references on the charging and discharging cooling powers $\dot{Q}_{TES}$ and $\dot{Q}_{TES,sec}$ will be quickly tracked, provided that the set points are achievable, which is ensured through the feasibility constraints given in \eqref{eq_Feasibility_Constraints}. Therefore, $\dot{Q}_{TES}$ and $\dot{Q}_{TES,sec}$ can be assumed to be constant during the whole sampling time $\Delta t$ and they match the reference values $\dot{Q}_{TES}^{ref}$ and $\dot{Q}_{TES,sec}^{ref}$ already applied to the system in the previous sampling time.

The thermal losses are not constant during the whole sampling time $\Delta t$, but they can be estimated from the temperature of the \emph{intermediate fluid}, for instance using a Tustin approximation, giving rise to energy losses $\Delta U_{surr}$. As a result, the energy balance can be expressed as:

\begin{equation}
	\begin{split}
		\Delta U_{TES}(t) &= \int_{0}^{\Delta t} \dot{Q}_{TES,int} \; dt \approx  \Delta \hat{U}_{TES}(t) = \\ 
		&= c_{p,TES,int} \, m_{TES,int} \, \Big[ T_{TES,int}(t) - T_{TES,int}(t-1) \Big] + \Big[\dot{Q}_{TES}^{ref}(t-1) - \dot{Q}_{TES,sec}^{ref}(t-1) \Big] \Delta t - \Delta U_{surr} (t) \;, \\
	\end{split}
	\label{eq_Energy_Estimator}
\end{equation}

\noindent where $\Delta \hat{U}_{TES}$ refers to the overall estimated energy transferred between the \emph{intermediate fluid} and all PCM cylinders during the complete sampling time $\Delta t$.

All the terms on the right-hand side of \eqref{eq_Energy_Estimator} are either known or can be computed when estimating the state vector $\mathbold{x}_{TES}(t)$. Once the overall transferred energy, $\Delta \hat{U}_{TES}(t)$, has been estimated, it is trivial to obtain the cold-energy transferred by every PCM cylinder $\Delta \hat{U}_{pcm}(t)$, and an algorithm identical to that described in subsection \ref{subSecPredictionModel} can be applied only for a sampling time ($k$ = 0), giving rise to the following estimation on the state vector,

\begin{equation}
	\mathbold{\hat{x}}_{TES}(t) = 
	\left[
		\begin{matrix} 
			\hat{h}_{pcm,1} (t) \\ 
			\hat{h}_{pcm,2} (t) \\ 
			\cdots \\ 
			\hat{h}_{pcm,n_{lay}} (t) \\ 
			T_{TES,int} (t)
		\end{matrix}
	\right] \;.
	\label{eq_IdentifiedStateVector}
\end{equation}

It is remarked that $T_{TES,int}(t)$ is not estimated but measured, and that is the key of the system feedback that allows to estimate the complete state vector based only on the simple energy balance on the \emph{intermediate fluid} shown in \eqref{eq_EnergyBalance}.

%%%%%%%%%%%%%%%%%%%%%%%%%%%%%%%%%%%%%%%%%%%%%%%%%%%%%%%%%%%%%%%%%%%%%%%%%%%%%%%%%%%%%

\section{Case study} \label{secCaseStudy}

In this section a case study is analysed in simulation, where a cooling demand profile that requires the combined power contribution of both the evaporator and the TES tank to face the peak demand is imposed. This case study turns out to be challenging, since the TES charging and discharging periods must be carefully scheduled to ensure the peak demand satisfaction, considering the plant power limits and the fact that the TES tank must remain within the latent zone. The actual energy costs throughout a certain day are also considered in the economic objective function. Some simulation results are presented and discussed, while the main advantages of the proposed scheduler with respect to the previous strategies are remarked. 

The design parameters of the TES tank are detailed in Table \ref{tabTESTankParameters}, while the most relevant thermodynamic properties of the PCM are specified in Table \ref{tabPCMProperties}. The intermediate fluid is a 60\% in volume ethylene glycol aqueous solution with very high thermal conductivity, while the secondary fluid is a 60\% in volume propylene glycol aqueous solution. Regarding the refrigeration cycle, it works with R404a as refrigerant. The models of the components described in the work by Bejarano \cite{bejarano2017optimization} are applied, whereas the steady-state parameters experimentally identified in the work by Bejarano \emph{et al.} \cite{bejarano2016identifying} have been used. Specifically, the parameters related to the main compressor, the air condenser, the evaporator related to the refrigerated chamber at -20\degree C, and the corresponding expansion valve \emph{EEV2} are applied, while the TES expansion valve is assumed to be identical to \emph{EEV2}. The thermodynamic properties of all fluids are computed using the \emph{CoolProp} tool \cite{CoolProp}.

\begin{table}[h]
	\centering
	\caption{Design parameters of the TES tank}
	\label{tabTESTankParameters}
	\begin{tabular} { L{1.5cm} L{7.5cm} L{2.0cm}  L{2.0cm} }
		\toprule
		\emph{Symbol} & \emph{Description} & \emph{Value} & \emph{Units} \\ 
		\toprule
		$L_{tank}$ & Length of the TES tank & 1.4 & m \\ 
		\midrule
		$D_{tank}$ & Internal diameter of the TES tank & 0.4 & m \\ 
		\midrule
		$e_{tank}$ & Thickness of the TES tank wall & 0.005 & m \\ 
		\midrule
		$n_{pcm}$ & Number of PCM cylinders & 17 & -- \\ 
		\midrule
		$D_{pcm}$ & External diameter of the PCM cylinders & 0.0445 & m \\ 
		\midrule
		$e_{pcm}$ & Thickness of the PCM cylinder coating & 0.001 & m \\ 
		\midrule
		$\kappa_{coat,pcm}$ & Thermal conductivity of the PCM cylinder coating & 16.3 & W m\textsuperscript{-1} K\textsuperscript{-1} \\ 
		\midrule
		$n_{refr}$ & Number of refrigerant pipes & 36 & -- \\ 
		\midrule
		$D_{refr}$ & External diameter of the refrigerant pipes & 0.020 & m \\ 
		\midrule
		$e_{refr}$ & Thickness of the refrigerant pipe wall & 0.001 & m \\ 
		\midrule
		$\kappa_{coat,refr}$ & Thermal conductivity of the refrigerant pipe wall & 16.3 & W m\textsuperscript{-1} K\textsuperscript{-1} \\ 
		\midrule
		$n_{sec}$ & Number of secondary fluid pipes & 32 & -- \\ 
		\midrule
		$D_{sec}$ & External diameter of the secondary fluid pipes & 0.020 & m \\ 
		\midrule
		$e_{sec}$ & Thickness of the secondary fluid pipe wall & 0.001 & m \\ 
		\midrule
		$\kappa_{coat,sec}$ & Thermal conductivity of the secondary fluid pipe wall & 16.3 & W m\textsuperscript{-1} K\textsuperscript{-1} \\ 
		\midrule
		$V_{TES,int}$ & Volume of the \emph{intermediate fluid} & 0.109 & m\textsuperscript{3} \\
		\midrule
		$\alpha_{surr}$ & Coefficient of thermal losses & 0.1 & W m\textsuperscript{-2} K\textsuperscript{-1} \\
		\bottomrule
	\end{tabular}
\end{table}

\begin{table}[h]
	\centering
	\caption{Phase change material properties}
	\label{tabPCMProperties}
	\begin{tabular} { L{1.5cm} L{6.5cm} L{2.0cm}  L{2.0cm} }
		\toprule
		\emph{Symbol} & \emph{Description} & \emph{Value} & \emph{Units} \\ 
		\toprule
		${c_{p}}_{pcm}$ & Specific heat at constant pressure & 3690 & J kg\textsuperscript{-1} K\textsuperscript{-1} \\ 
		\midrule
		$h_{pcm}^{lat}$ & Specific enthalpy of fusion (latent phase) & 222000 & J kg\textsuperscript{-1} \\ 
		\midrule
		$T_{pcm}^{lat}$ & Phase change temperature & -29 & \degree C \\ 
		\midrule
		$\kappa_{pcm}$ & Thermal conductivity & 0.64 & W m\textsuperscript{-1} K\textsuperscript{-1} \\ 
		\midrule
		$\rho_{pcm}$ & Density & 1420 & kg m\textsuperscript{-3}\\ 
		\bottomrule
	\end{tabular}
\end{table}

\subsection{Cooling demand profile} \label{subSecCoolingDemandProfile}

\begin{figure}[h]
	\centering
	\includegraphics[width=7.0cm] {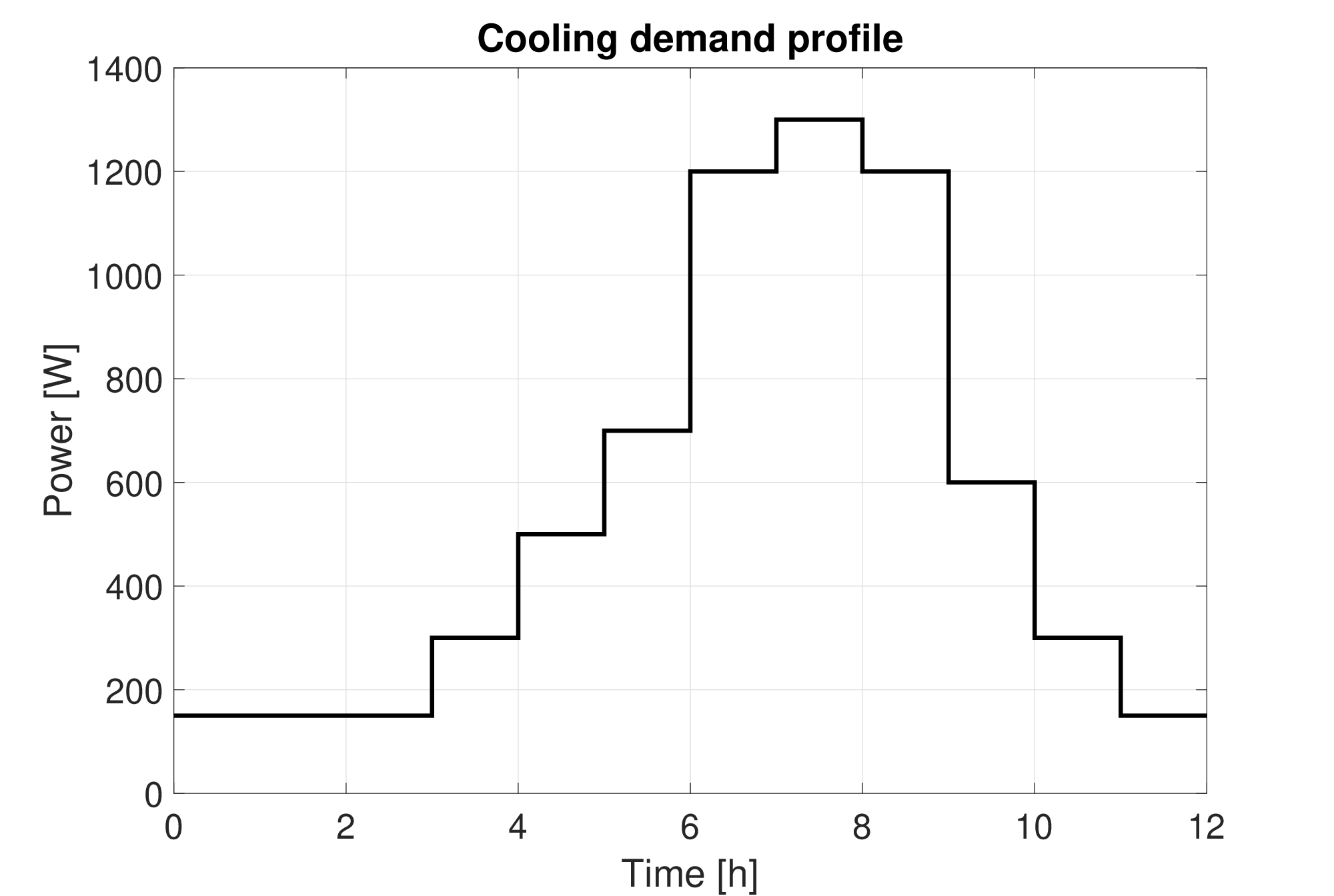}
	\caption{Cooling demand profile.} 
	\label{figCoolingDemandProfile}
\end{figure}

The realistic daily cooling demand profile represented in Figure \ref{figCoolingDemandProfile} is analysed hereafter. Regarding the peak demand, it has been tailored to the maximum combined cooling power achievable by the system, already considered in a previous work \cite{bejarano2019NMPC}. Moreover, the time window has been reduced to 12 hours instead of a complete day, according to the maximum charging and discharging periods considered in the design stage of the TES tank. The latter was described in detail in a previous work by Bejarano \emph{et al.} \cite{Bejarano2017NovelSchemePCM}, where the only difference was related to the PCM encapsulation. Since it is a research facility, 3-4 hour periods for full charging/discharging were regarded in the design stage as most desirable. However, since only a time scaling has been applied, the conclusions drawn are applicable to 24-hour operation, provided that the TES tank is designed accordingly.

It is observed in Figure \ref{figCoolingDemandProfile} that the cooling demand is not zero at any instant throughout the day, as in the case of industrial refrigeration or supermarket operation. A minimum value is required even during the \emph{night} hours, whereas the peak demand happens one hour past noon ($t$ = 7 h). This peak demand requires the combined contribution of the cooling power provided at the evaporator and the TES tank discharge, and thus it is expected that the TES tank is charged during the \emph{night} hours to store enough cold energy to face the discharging process surely imposed during the \emph{daylight} hours around noon. Anyway, as stated in Sections \ref{secIntroduction} and \ref{secScheduling}, the operating mode scheduling is also included in the optimization and it will be optimally set to satisfy the cooling demand while observing the imposed constraints and minimizing the operating cost.

\subsection{Energy price} \label{subSecEnergyPrice}

\begin{figure}[h]
	\centering
	\includegraphics[width=7.0cm] {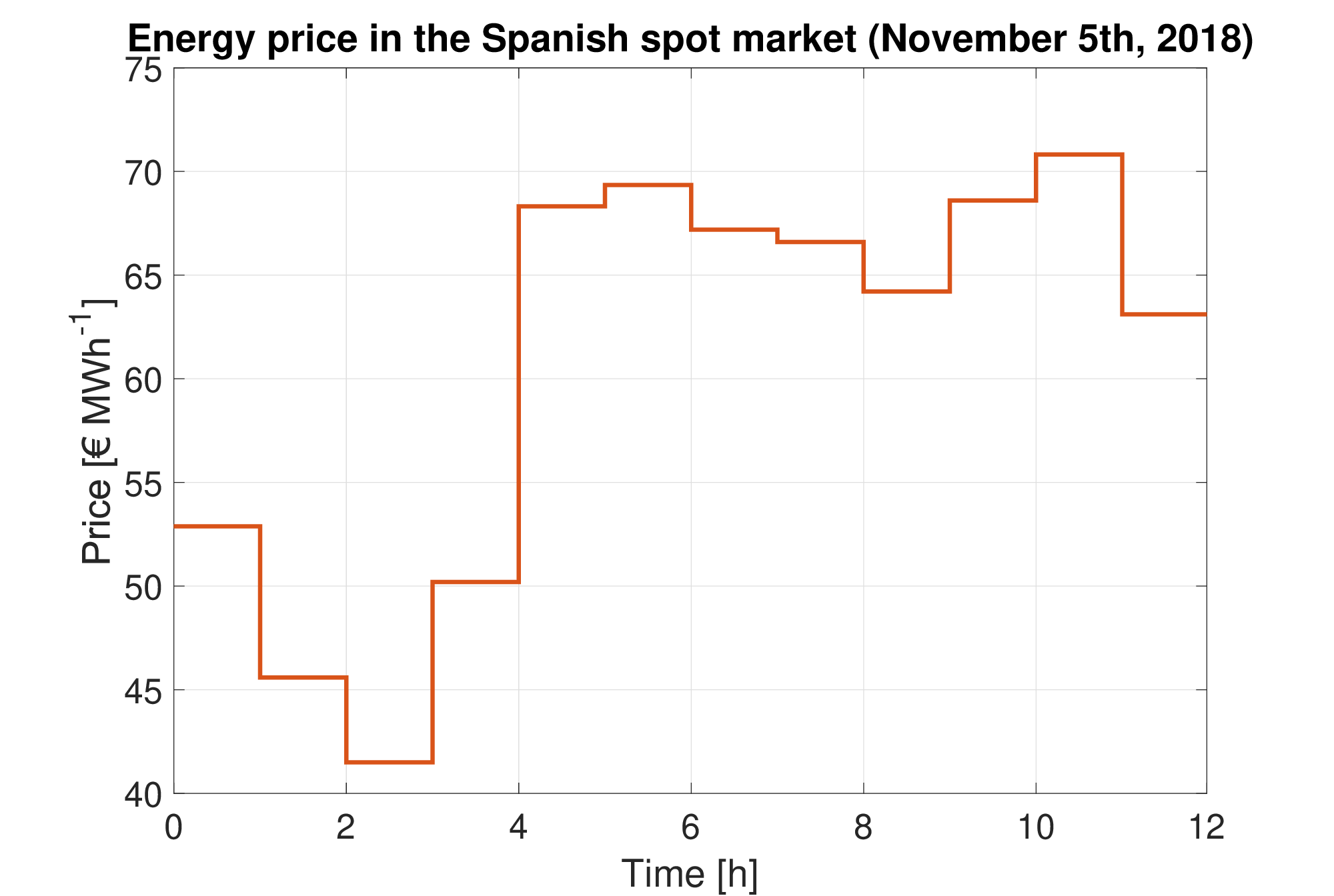}
	\caption{Energy price profile corresponding to November 5th, 2018 in the Spanish spot market \cite{REN_Portugal}.} 
	\label{figEnergyPrices}
\end{figure}

Actual energy prices are also considered in the objective function of the optimization problem, shown in \eqref{eq_Objective_Function}. Figure \ref{figEnergyPrices} shows the energy prices corresponding to a given day (November 5th, 2018) in the Spanish spot market, once again time-scaled into a 12-hour timeframe \cite{REN_Portugal}.

These prices, suitably scaled, have been applied to the cooling power generation at the evaporator, $w_{e,sec}$, and to the TES charging cooling power, $w_{TES}$. However, the cost related to the TES discharging cooling power $w_{TES,sec}$ has been set to zero, since releasing the previously stored cold energy does not involve a direct economic cost, always without considering the electrical power devoted to impulsing the secondary fluid.

\subsection{Simulation results} \label{subSecSimulationResults}

Some simulation results of the proposed scheduling strategy are presented in this subsection. A sampling time of 1 h has been selected, according to the cooling demand profile shown in Figure \ref{figCoolingDemandProfile}. Safety limits $\gamma_{TES}^{min}$ = 0.05 and $\gamma_{TES}^{max}$ = 0.95 have been imposed, while a prediction horizon of 12 h has been considered in the MINLP-based scheduling strategy. This choice is motivated by the fact that if the energy price and the cooling demand forecasts are given for a complete day, then the prediction horizon should cover at least this period. In the (rare) case that those forecasts vary widely between two consecutive days, it would be useful to consider a wider prediction/control horizon, but the computational load of the proposed strategy would be compromised, since the number of decision variables is shown to be proportional to the prediction/control horizon. Regarding the optimization tool, the OPTI Toolbox \cite{OPTI_ToolBox} has been applied in the MATLAB\textsuperscript{\textregistered} environment, while the BONMIN algorithm has been used to solve the mixed integer non-linear program \cite{bonami2008algorithmic}. BONMIN uses the Interior Point OPTimizer (IPOPT) for solving relaxed problems and Coin-OR Branch and Cut (CBC) as the mixed integer solver \cite{wachter2006implementation,CBC_Algorithm}.

\begin{figure}[h]
	\centering
	\includegraphics[width=7.0cm] {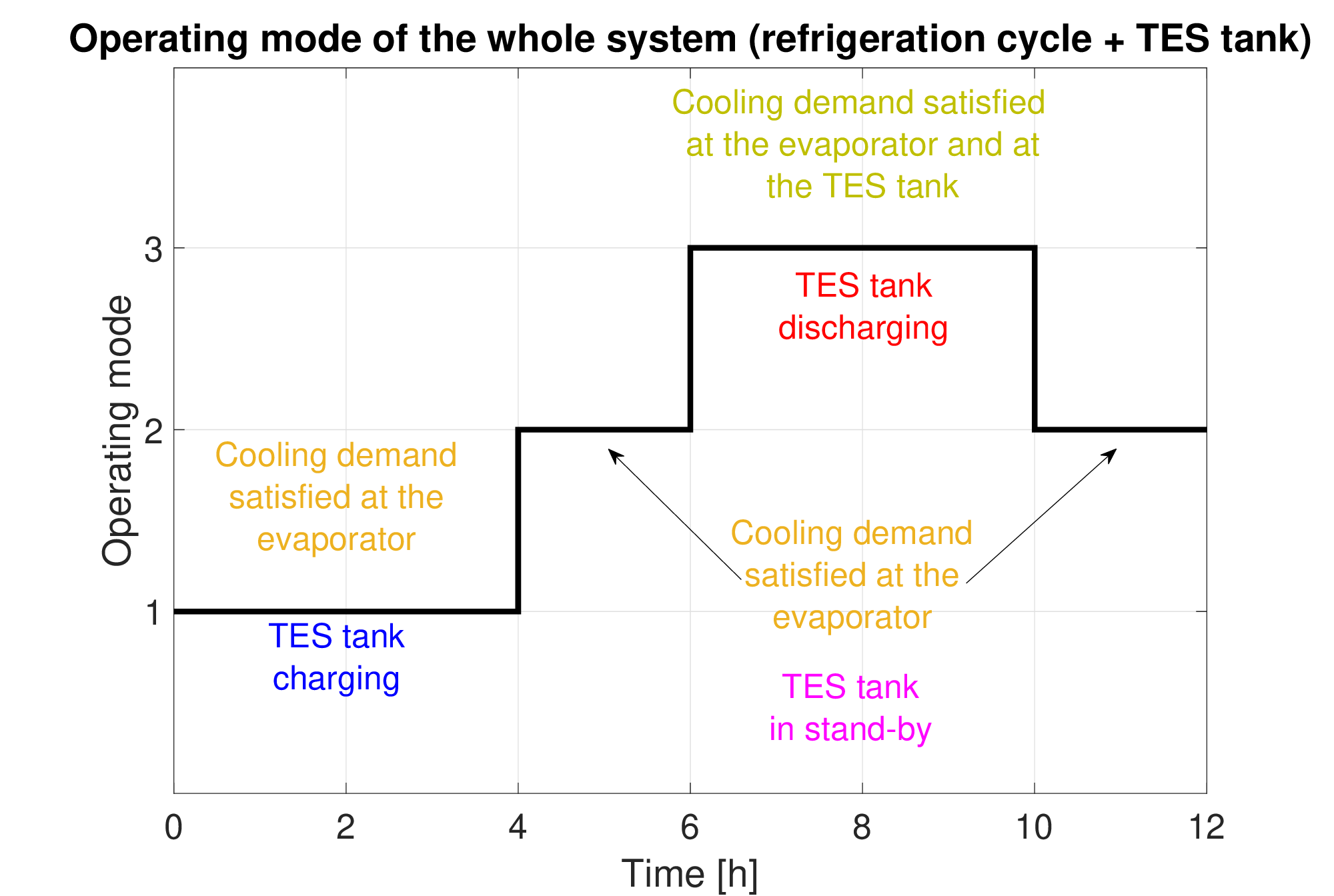}
	\caption{Optimal operating mode scheduling of the whole system given by the MINLP-based strategy.}	
	\label{figOperatingMode}				
\end{figure}

The optimal operating mode scheduling given by the mixed-integer non-linear optimization is represented in Figure \ref{figOperatingMode}, computed from the binary variable set \{$\delta_{e,sec}$, $\delta_{TES}$, $\delta_{TES,sec}$\} defining whether the corresponding cooling power is active or not, according to subsection \ref{subSecOperatingModes}.

It is confirmed in Figure \ref{figOperatingMode} that the intuitive operating mode scheduling suggested in subsection \ref{subSecCoolingDemandProfile} is actually set by the MINLP-based strategy as the optimal one, given the actual energy price variations described in subsection \ref{subSecEnergyPrice}. However, two different periods where the cooling demand is satisfied only at the evaporator are also included, separating the TES tank charging and discharging processes. The constraints regarding the operating mode limitation described in \eqref{eq_BinaryConstraints} are shown to be observed, since only modes 1 to 3 have been scheduled.

\begin{figure}[h]
	\centering
	\subfigure[Reference on $\dot{Q}_{e,sec}$.]{
		\includegraphics[width=6.5cm] {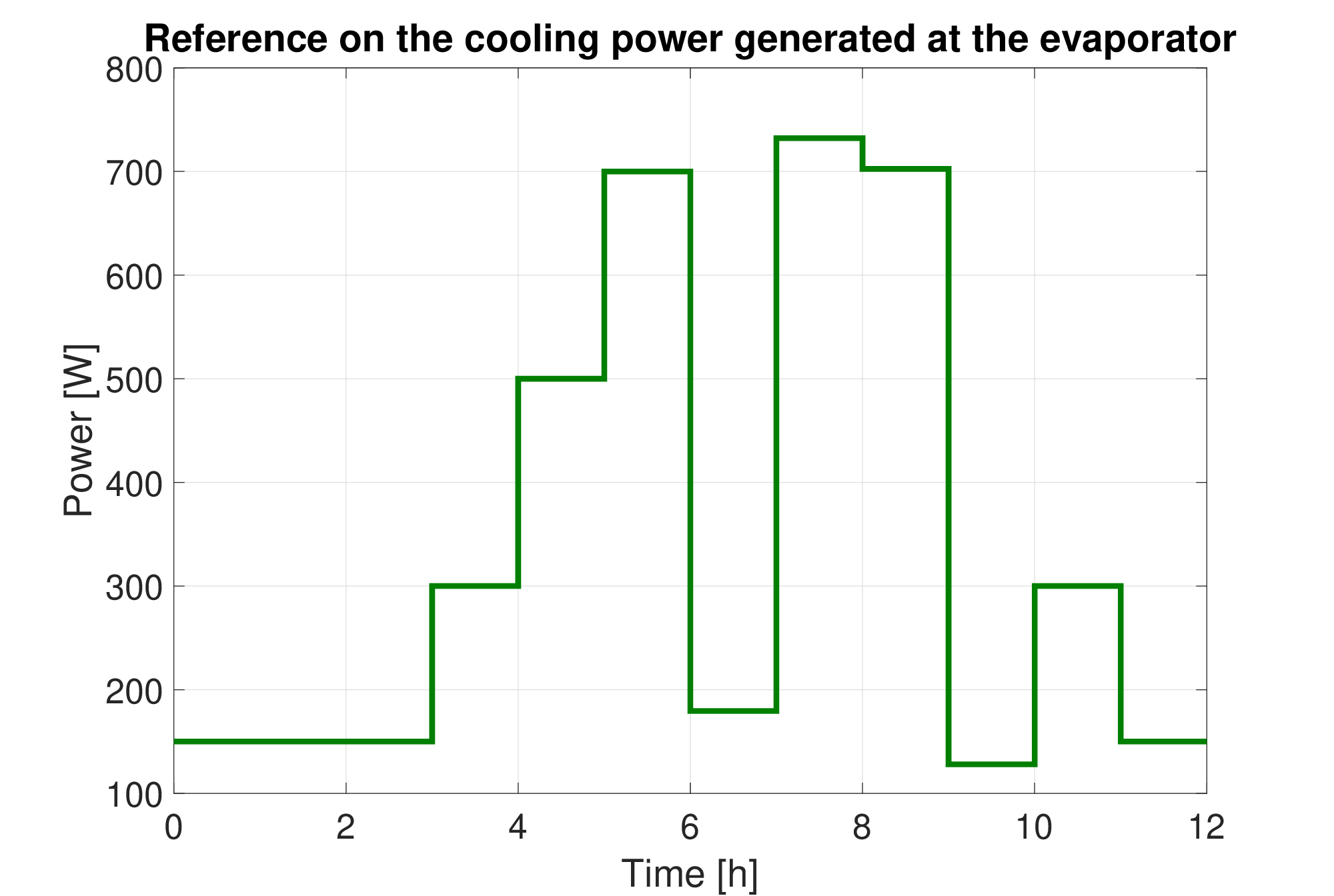}
		\label{figQ_e_sec}
	}
	\subfigure[Reference on $\dot{Q}_{TES}$.]{
		\includegraphics[width=6.5cm] {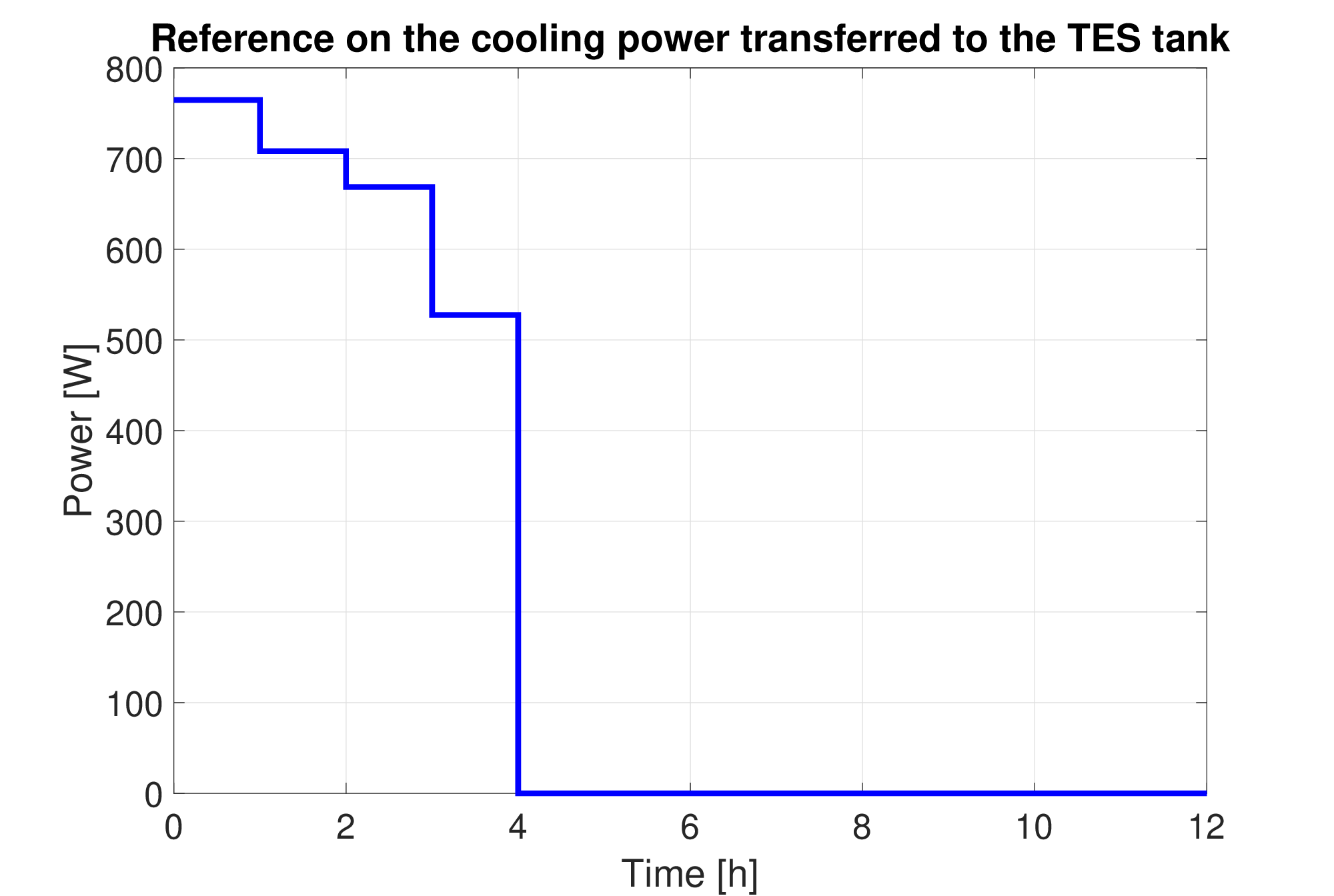}
		\label{figQ_PCM}
	}
	\subfigure[Reference on $\dot{Q}_{TES,sec}$.]{
		\includegraphics[width=6.5cm] {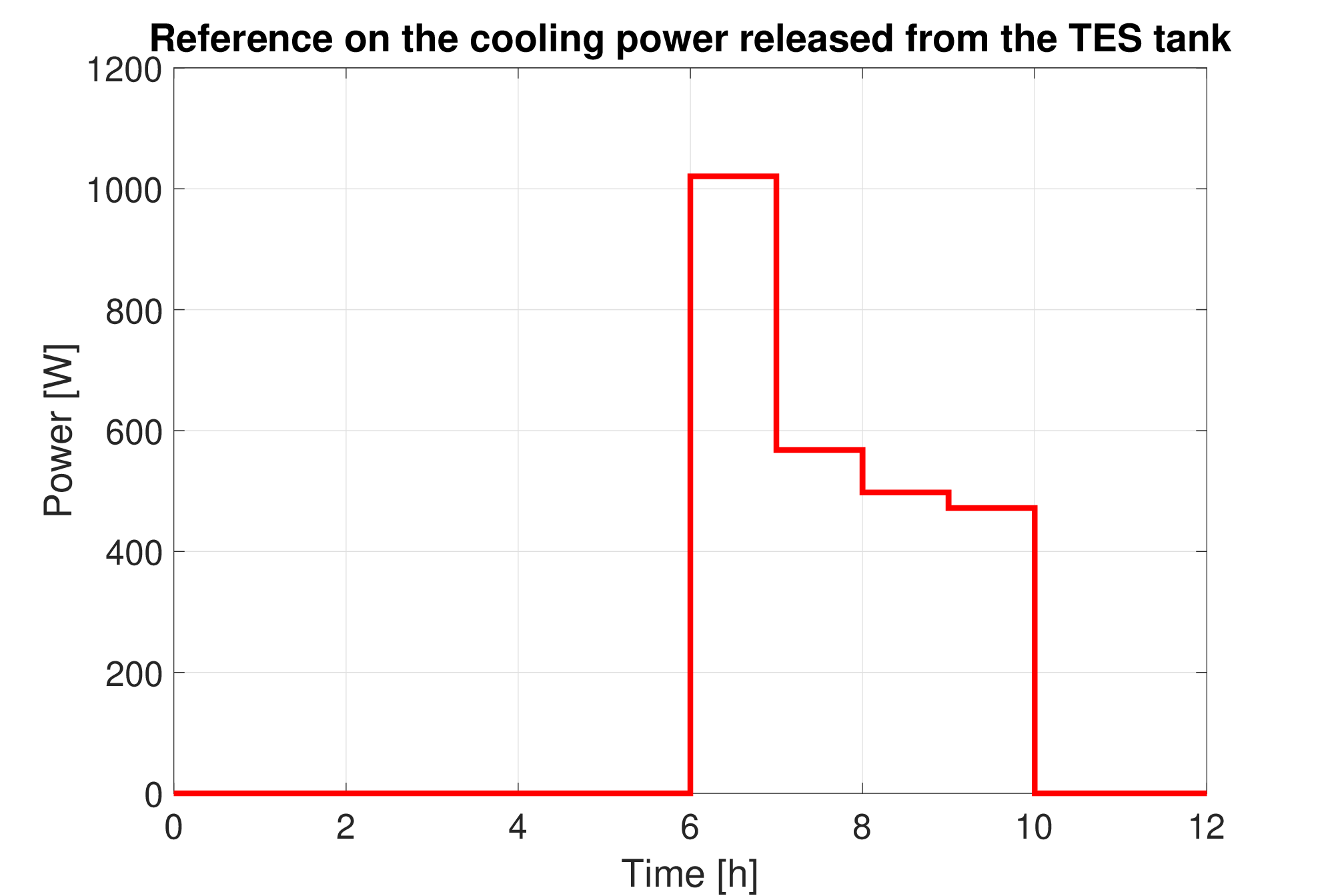}
		\label{figQ_PCM_sec}
	}	
	\caption{References on the relevant cooling powers given by the proposed MINLP-based scheduling strategy.}
	\label{figCoolingPowers}				
\end{figure}

\begin{figure}[h]
	\centering
	\includegraphics[width=6.5cm] {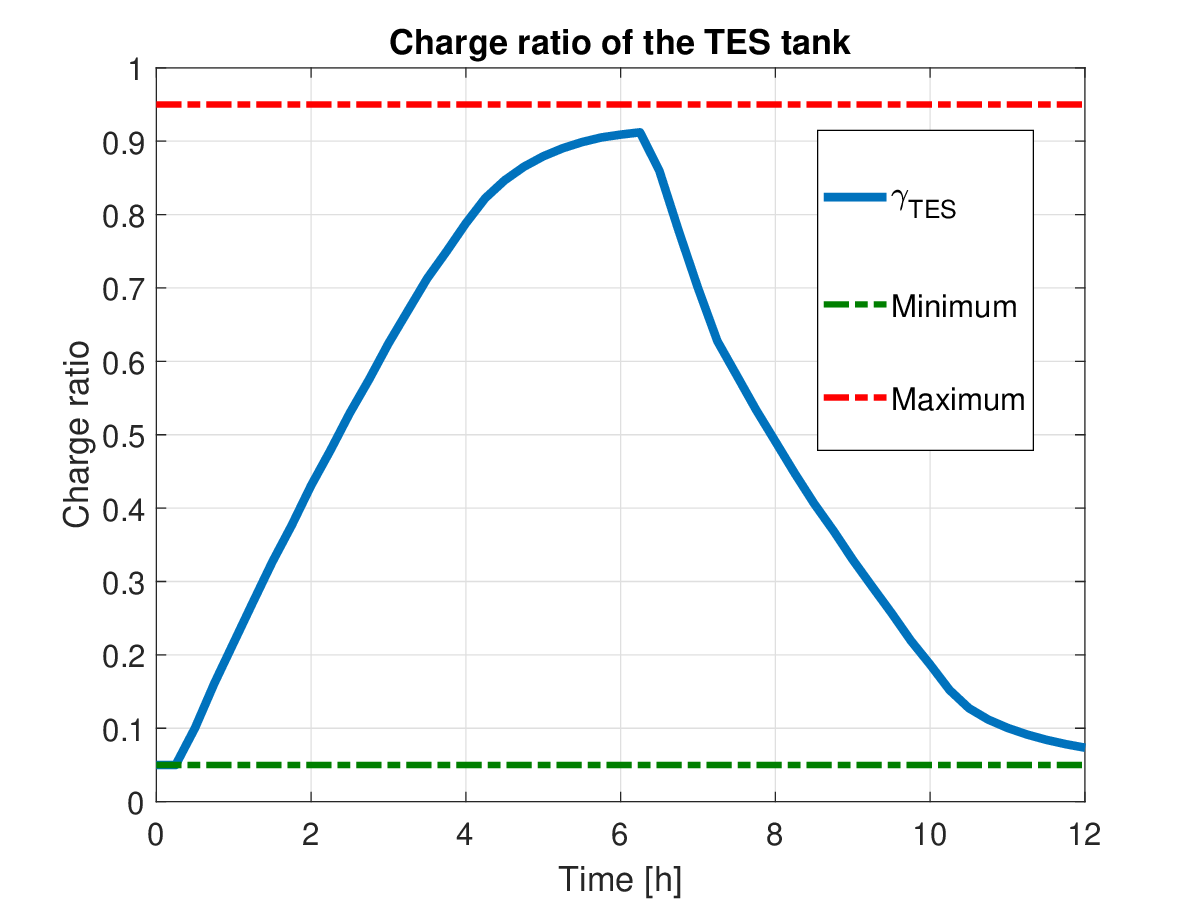}
	\caption{\emph{Charge ratio} of the TES tank for the proposed MINLP-based scheduling strategy.} 
	\label{figChargeRatio}
\end{figure}

\begin{figure}[h]
	\centering
	\includegraphics[width=9.5cm] {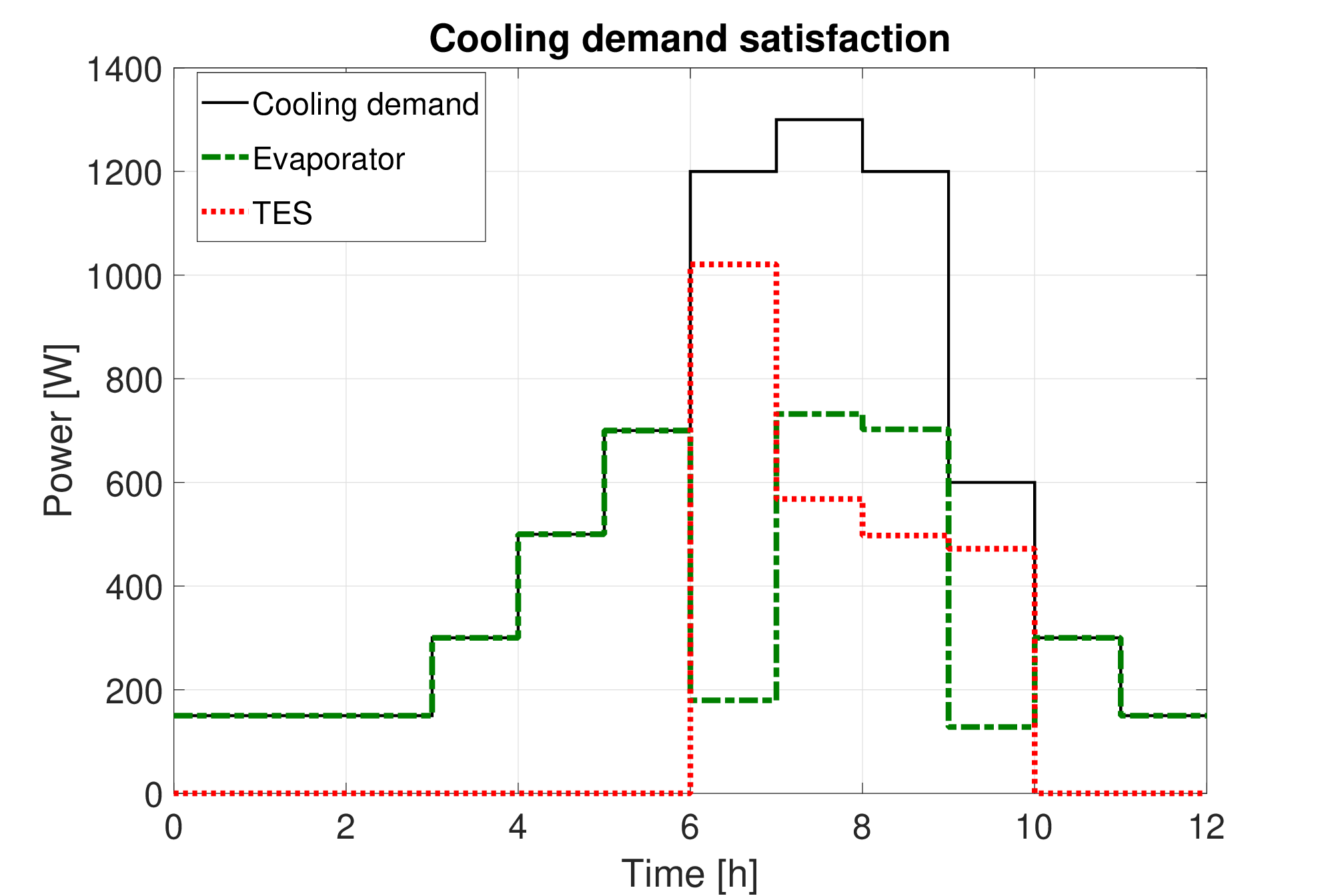}
	\caption{Cooling demand satisfaction of the proposed MINLP-based scheduling strategy.}	
	\label{figCoolingDemandSatisfaction}				
\end{figure}

The references on $\dot{Q}_{e,sec}$, $\dot{Q}_{TES}$, and $\dot{Q}_{TES,sec}$ are represented in Figure \ref{figCoolingPowers}, while the evolution of the TES tank \emph{charge ratio} throughout the day is shown in Figure \ref{figChargeRatio}. Eventually, the cooling demand satisfaction is shown in Figure \ref{figCoolingDemandSatisfaction} by representing together both contributions to the cooling power provided to the secondary fluid: $\dot{Q}_{e,sec}$ and $\dot{Q}_{TES,sec}$.

The optimal references on the relevant cooling powers shown in Figure \ref{figCoolingPowers} allow to ensure the satisfaction of the cooling demand profile represented in Figure \ref{figCoolingDemandProfile} all throughout the day, as shown in Figure \ref{figCoolingDemandSatisfaction}, as a result of the application of \eqref{eq_Cooling_Demand_Constraint}. Moreover, the references  $\dot{Q}_{e,sec}^{ref}$, $\dot{Q}_{TES}^{ref}$, and $\dot{Q}_{TES,sec}^{ref}$ meet the power limits imposed through the constraints indicated in \eqref{eq_Feasibility_Constraints}, and thus they are achievable by the TES-backed-up refrigeration cycle. Eventually, the TES tank remains within the security range given by $\gamma_{TES}^{min}$ and $\gamma_{TES}^{max}$, as may be checked in Figure \ref{figChargeRatio}, thus complying with the security limits indicated in \eqref{eq_GammaTESConstraints}.

The main advantages of the proposed strategy with respect to the NMPC-based scheduling are related to three items:

\begin{enumerate}[(a)]
	
	\item Operating cost: the daily operating cost is minimized, as shown in \eqref{eq_Objective_Function}. Indeed, the operating mode scheduling given by the MINLP-based strategy is intended to be optimal, and it might not be completely intuitive as previously discussed in the view of the cooling demand profile shown in Figure \ref{figCoolingDemandProfile}. In the NMPC-based scheduling strategy, it is necessary to suggest a certain operating mode scheduling \emph{a priori}, given the cooling demand profile and energy prices. The suggested scheduling might not be optimal, which surely involves a higher daily operating cost.
	
	\item Ease of tuning: it was stated in the work by Bejarano \emph{et al.} \cite{bejarano2019NMPC} that some weights in the objective function related to the \emph{charge ratio} must be tuned in order to promote the TES tank charging/discharging, according to the suggested operating mode scheduling. These weights may be difficult to tune in some cases and it is easy to perceive that different tuning may lead to problem infeasibility if the cooling profile is demanding enough and the TES tank is not properly charged/discharged when the corresponding operating modes are scheduled, that would result in cooling demand non-satisfaction.  
	
	\item Adaptability to cooling demand variations: in the NMPC-based strategy, the operating mode scheduling must be proposed in the light of the forecast on the cooling demand and energy prices. Therefore, if the predicted cooling demand changes, the suggested operating mode scheduling is likely not to be optimal, at best, but at worst it may also lead to problem infeasibility, that would also result in cooling demand non-satisfaction.  
	
\end{enumerate}

If the optimal operating mode scheduling represented in Figure \ref{figOperatingMode} is imposed in the NMPC-based strategy, and the weights in the objective function related to the \emph{charge ratio} are carefully tuned, the optimal daily operating cost and the same cooling power references as those shown in Figure \ref{figCoolingPowers} can be achieved, but the difficulty in getting this tuning is very high when compared to the reduced set of tuning parameters of the MINLP-based strategy. 

%%%%%%%%%%%%%%%%%%%%%%%%%%%%%%%%%%%%%%%%%%%%%%%%%%%%%%%%%%%%%%%%%%%%%%%%%%%%%%%%%%%%%

\section{Conclusions and future work} \label{secConclusions}

In this brief, the operation of a hybrid system consisting of a vapour-compression refrigeration cycle and a PCM-based TES unit has been analysed. The work has been focused on the scheduling problem arising when a certain demand profile is imposed and the references on the cooling powers involved (TES charging/discharging and power provided at the evaporator) are intended to be optimally scheduled, according to energy price forecast, feasibility constraints, and limited storable cold energy. The application of the different operating modes to this problem has been discussed, whereas a subset has been considered as most likely to be scheduled. 

The proposed scheduling strategy based on the predictive control paradigm has been described in detail. Concerning the prediction model, a previously presented simplified model focused on the dominant dynamics related to heat transfer within the TES tank has been used as a starting point, and a first-order linearisation based on the ideas of the PNMPC has been applied to compute a even more simplified model, suitable to be used within the optimization algorithm. A hybrid decision set including both binary and continuous variables is considered, whereas several constraints concerning cooling demand satisfaction, power feasibility, operating mode limitation, and TES latency are applied. Since the decision set include binary and continuous variables, mixed-integer non-linear programming is needed to solve the optimization problem, where the economic operating cost is minimized. State estimation has also been addressed, since the TES state vector is not fully measurable. An energy balance on the TES \emph{intermediate fluid} has been applied to estimate the energy transferred between the latter and the PCM, while a step-by-step procedure has been proposed to estimate the complete TES state vector, including the enthalpy distribution within the PCM cylinders.

Some simulation results have been presented for a challenging cooling demand profile that forces the optimizer to schedule alternative TES charging and discharging processes to face the peak demand. The proposed strategy is shown to provide the optimal operating mode scheduling as well as the references on the cooling powers involved, that meet the feasibility constraints and satisfy the cooling demand. Moreover, the TES \emph{charge ratio} is shown to remain within the latency limits imposed. The MINLP-based strategy is shown to be much easier to tune than the NMPC-based strategy. Furthermore, the first one is shown to adapt to demand variations. Both advantages allow to ensure the cooling demand satisfaction and problem feasibility, while the daily operating cost is ensured to be minimized.

As future work, the proposed scheduling and control strategy is planned to be applied to the experimental facility as soon as it is fully operative. Furthermore, the variations of the system performance (\emph{e.g.} COP) in the different operating modes should be considered in the objective function.

\section*{Acknowledgement}

\ack{The authors would like to acknowledge Spanish MCeI (Grants DPI2015-70973-R and DPI2016-79444-R) for funding this work, as well as University of Seville through VI PPIT-US program. The cooperation of INESC-ID was supported by FCT (Portugal) under UID/CEC/50021/2019, and by POR Lisboa-Lisboa-01-0145-FEDER-031411.}

\bibliography{bibliografia}

\end{document}